\documentclass[11pt,reqno]{amsart}
\usepackage{amsaddr}
\usepackage{graphicx}
\usepackage{amssymb}
\usepackage{amsfonts}
\usepackage[]{amsmath}
\usepackage[]{epsfig}
\usepackage[]{pstricks}
\newpsobject{malla}{psgrid}{subgriddiv=1,griddots=10,gridlabels=6pt}
\usepackage[]{float}
\usepackage{setspace}
\usepackage{tikz}

\newpsobject{malla}{psgrid}{subgriddiv=1,griddots=10,gridlabels=6pt}
\newtheorem{theorem}{Theorem}[section]
\newtheorem{lemma}{Lemma}[section]
\newtheorem{definition}{Definition}[section]

\newtheorem{corollary}{Corollary}[section]

\numberwithin{equation}{section}

\theoremstyle{remark}
\newtheorem{remark}{Remark}[section]
\newcommand{\R}{\mathbb R}

\begin{document}

	\pagenumbering{arabic}	
\title[KdV on a metric star graph]{The Korteweg-de Vries equation on a metric star graph}

\author{M\'arcio CAVALCANTE}
\address{\emph{Instituto de Matem\'atica, Universidade Federal de Alagoas,\\ Macei\'o-Brazil}}
\email{marcio.melo@im.ufal.br}
\let\thefootnote\relax
\footnote{Permanent address: Instituto de Matem\'atica, Universidade Federal de Alagoas,
	Macei\'o-Brazil
}

\begin{abstract}
We prove local well-posedness for the Cauchy problem associated
to Korteweg-de Vries equation on a metric star graph  with three semi-infinite edges given by one negative half-line and two positives half-lines attached to a common vertex, for two classes of boundary conditions. The results are obtained
in the low regularity setting by using the Duhamel Boundary Forcing Operator, in context of half-lines, introduced by Colliander, Kenig (2002), and extended by Holmer (2006) and Cavalcante (2017).

\end{abstract}
\maketitle

\section{Introduction}
Partial differential equations (PDE) on metric graphs, known as networks or quantum
graphs, arise naturally in many topics of physics such as acoustics, optics, condensed
matter and polymer physics, and,
more recently, in connection with biological trees. We refer to \cite{Berkolaiko} and \cite{Mehmeti} for further information and bibliography. Earlier, the
linear Schr\"odinger equation on a metric graph was subject of extensive research due to its applications in quantum  nanotechnologies, chemistry and mesoscopic
physics (see \cite{Bellazi}, \cite{Berkolaiko} and references therein). Studies of the nonlinear Schr\"odinger
equation on graphs have started appearing recently. In particular, existence and
stability of standing waves for nonlinear Schr\"odinger equation on a star graph with a power nonlinearity $|u|^{p-1}u$ have been studied extensively, see \cite{Adami}, \cite{Cacciapuoti} and  \cite{Noja} for a brief survey of the topic.  Recently, Ardila \cite{Ardila} has obtained existence and stability of standing waves for the logarithmic Schr\"odinger equation on star graphs. Caudrelier \cite{Caudrelier} presented a method to solve the open problem of formulating the inverse scattering method for an integrable PDE on a star graph and the nonlinear Schr\"odinger equation was chosen
to illustrate the method. 

 Another nonlinear dispersive equation, the Benjamin-Bona-Mahony (BBM) equation, is treated in \cite{Bona} and  \cite{Mugnolo}. More precisely,  Bona and Cascaval \cite{Bona} obtained local well-posedness in Sobolev space $H^1$ and  Mugnolo and Rault \cite{Mugnolo} obtained existence of traveling waves for the BBM equation on star graphs.

The well-known Korteweg-de Vries (KdV) equation
\begin{equation}\label{KdV}
u_t+u_{xxx}+u_xu=0
\end{equation}
was first derived by Korteweg and de Vries \cite{KDV} in 1895 as a model for long waves propagating on a
shallow water surface.  It is now commonly accepted as a mathematical model
for the unidirectional propagation of small-amplitude long waves in nonlinear dispersive
systems. In particular, the KdV equation is not only used to serve as a model to study surface
water waves. In fact, recently the KdV equation has been used to
serve as a model to study blood pressure waves in large arteries. In this way, for example, in \cite{Chuiko}  a new computer model for systolic pulse waves within the cardiovascular system based on the  KdV equation is presented. In \cite{Crepeau} some particular solutions of the KdV equation, the $2$ and $3$-\textit{soliton} solutions, seem to be good candidates to match the observed pressure pulse waves. 

In the mathematic context the Cauchy problem for the KdV posed on the
real axis, torus, on the half-lines and on a finite interval have
 been well studied in the last years, we refer as a sample \cite{BSZ2, Bonaint, CK, Faminskii, Guo, Holmer, Jia, KPV, Kishimoto}, but we point out that there are many further references in the works just cited. We finally notice that Cavalcante and Mu\~noz \cite{CM} showed
 that \textit{solitons} posed initially far away from the origin are strongly stable for
 the IBVP associated to the KdV equation posed on the right half-line, assuming homogeneous boundary
 conditions.

\subsection{Previous works for the linearized KdV equation on star graphs} Studies for the linearized Korteweg-de Vries equation on star graphs have started appearing recently. In \cite{Sobirov1},\cite{Sobirov2} and \cite{Sobirov3} was studied existence and uniqueness of solutions for the linearized KdV equation on metric star graphs by using potential theory, where the solutions were obtained in the class of Schwartz and in Sobolev classes with high order. Very recently, Mugnolo, Noja and Seifert \cite{Noja2} obtained a characterization of all boundary conditions  under which the Airy-type evolution equation $u_t=\alpha u_{xxx}+\beta u_x$, for $\alpha\in \R\setminus \{0\}$ and $\beta\in \R$ on star graphs given by $E=E_{-}\cup E_{+}$, where $E_{+}$ and $E_{-}$ are finite or countable collections of semi-infinite edges parametrized by $(-\infty,0)$ or $(0,+\infty)$ respectively and the half-lines are connected at a vertex $v$, generates a semigroup of contractions.

 As far as we know, the study of the nonlinear Korteweg-de Vries equation in star graphs is unknown. Here we extend the treatment
of the KdV equation given in \cite{Holmer} from intervals and half-lines to a star graph with three edges. 

\subsection{Formulation of the problem}In this work, motivated by the work \cite{Sobirov2}  we consider the KdV equation on a star graph $\mathcal{Y}=(-\infty,0)\cup (0,+\infty)\cup (0,+\infty)$ with three semi-infinite edges given by one negative half-line and two positives half-lines attached to a common vertex, also known as $\mathcal{Y}$-junction.  More precisely,
\begin{equation}\label{KDV}
\begin{cases}
u_t+u_{xxx}+u_xu=0,& (x,t)\in(-\infty,0)\times (0,T),\\
v_t+v_{xxx}+v_xv=0,& (x,t)\in(0,+\infty)\times (0,T),\\
w_t+w_{xxx}+w_xw=0,& (x,t)\in(0,+\infty)\times (0,T).
\end{cases}
\end{equation}
We consider the initial conditions given by
\begin{equation}\label{initial}
u(x,0)=u_0(x),\;v(x,0)=v_0(x)\; \text{and}\;w(x,0)=w_0(x),
\end{equation}
where
\begin{equation}\label{spacetrace}
(u_0,v_0,w_0)\in H^{s}(\R^-)\times H^{s}(\R^+)\times H^{s}(\R^+):=H^{s}(\mathcal{Y}).
\end{equation}
Our goal in studying Cauchy problem \eqref{KDV}-\eqref{initial} is to obtain results of local well-posedness, in  Sobolev spaces with low regularity.

\subsection{Choices of boundary conditions} Determining the number of boundary conditions necessary for a well-posed problem is a nontrivial
issue. As far we know, it's not at all clear which boundary conditions should be appropriate for physical applications, and therefore here we will consider two classes of boundary conditions that are coherent with uniqueness calculations for smooth decaying solutions of a linear version of the Cauchy problem \ref{KDV}. In this sense, suppose that $(u(x,t),v(x,t), w(x,t))$ is a smooth decaying solution   of a linear version of \eqref{KDV}, i.e.
\begin{equation}\label{KDVlinear}
	\begin{cases}
		u_t+u_{xxx}=0,& (x,t)\in(-\infty,0)\times (0,T),\\
		v_t+v_{xxx}=0,& (x,t)\in(0,+\infty)\times (0,T),\\
		w_t+w_{xxx}=0,& (x,t)\in(0,+\infty)\times (0,T),
	\end{cases}
\end{equation}
with homogeneous initial condition $(u_0,v_0,w_0)=(0,0,0).$  Multiplying the equations in \eqref{KDVlinear} by $u,$ $v$ and $w$ respectively, and integrating by parts we obtain
\begin{equation}\label{uniqueness}
\begin{split}
&\int_{-\infty}^0u^2(x,T)dx+\int_{0}^{+\infty}v^2(x,T)dx+\int_{0}^{+\infty}w^2(x,T)dx\\
&\quad \quad=\int_0^T (u_x^2(0,t)-v^2_x(0,t)-w_x^2(0,t))dt\\
&\quad \quad \quad-2\int_0^Tu(0,t)u_{xx}(0,t)dt  +2\int_0^Tv(0,t)v_{xx}(0,t)dt +2\int_0^Tw(0,t)w_{xx}(0,t)dt.
\end{split}
\end{equation}

By analyzing \eqref{uniqueness}, we are interested in boundary conditions for Cauchy problem \eqref{KDV}-\eqref{initial} such that the right hand side of \eqref{uniqueness} would have a non positive sign. 

In this sense, if we consider the following particular boundary conditions
\begin{equation}\label{v1}
u(0,t)=\alpha_2v(0,t)=\alpha_3w(0,t),\; t\in (0,T),
\end{equation}
\begin{equation}\label{v2}
u_x(0,t)=\beta_2v_x(0,t)+\beta_3w_x(0,t),\; t\in (0,T)
\end{equation}
and
\begin{equation}\label{v3}
u_{xx}(0,t)=\frac{1}{\alpha_2}v_{xx}(0,t)+\frac{1}{\alpha_3}w_{xx}(0,t),\; t\in (0,T),
\end{equation}
where $\alpha_2,\ \alpha_3,\ \beta_2$ and $\beta_3$ are real constants satisfying  $3\beta_i^2\leq 1$ for $i=2,3$ we have that

\begin{equation}
\begin{split}
&\int_{-\infty}^0u^2(x,T)dx+\int_{0}^{+\infty}v^2(x,T)dx+\int_{0}^{+\infty}w^2(x,T)dx\\
&\quad=\int_0^T[(\beta_2v_x(0,t)+\beta_3w_x(0,t))^2-v^2_x(0,t)-w_x^2(0,t)]dt\\
&\quad \quad \quad+2\int_0^Tu(0,t)(-u_{xx}(0,t)+\frac{1}{\alpha_2}v_{xx}(0,t) +\frac{1}{\alpha_3}w_{xx}(0,t))\\
&\quad= \int_{0}^T[\beta_2^2v_x^2(0,t)+2\beta_2\beta_3v_x(0,t)w_x(0,t)+\beta_3^2w_x^2(0,t)-v^2_x(0,t)-w_x^2(0,t)]dt\\
&\quad \leq\int_0^T [\beta_2^2v_x^2(0,t)+2\beta_2^2v_x^2(0,t)+2\beta_3^2w_x^2(0,t)+\beta_3^2w_x^2(0,t)-v^2_x(0,t)-w_x^2(0,t)]dt\\
&\quad \leq\int_0^T [(3\beta_2^2-1)v_x^2(0,t)+(3\beta_3^2-1)w_x^2(0,t)]dt,
\end{split}
\end{equation}
where we have used the Cauchy-Schwarz inequality. It follows that $u(x,T)=v(x,T)=w(x,T)= 0$, which implies the uniqueness argument.

\begin{remark}
	The boundary conditions \eqref{v1}-\eqref{v3} also were considered by Sobirov, Akhmedov and Uecker \cite{Sobirov2} to treat a linearized version of KdV on a $\mathcal{Y}$-junction. Furthermore, this boundary condition is an example included in the general boundary conditions obtained  by Mugnolo, Noja and Seifert on the study of Airy equation on star graphs (see Example 4.10 in \cite{Noja2}). 
\end{remark}
\begin{remark}
	The calculation above is only formal, however in the context of half-lines this is very efficient to describe the boundary conditions for initial boundary value problems (IBVPs). In this sense, Holmer \cite{Holmer} describes boundary conditions for IBVPs associated to the KdV equation on the positive and negative half-lines that imply a result of  well-posedness. For other nice discussion about the boundary condition for KdV on half-lines, based on the behavior of characteristic curves we refer the reader to \cite{Deconinck}.
\end{remark}
In the same way, the following particular boundary conditions
\begin{equation}\label{2v1}
u(0,t)=\frac{1}{\alpha_2}v(0,t)+\frac{1}{\alpha_3}w(0,t),\; t\in (0,T),
\end{equation}
\begin{equation}\label{2v2}
u_x(0,t)=\beta_2v_x(0,t)+\beta_3w_x(0,t),\; t\in (0,T)
\end{equation}
and
\begin{equation}\label{2v3}
u_{xx}(0,t)=\alpha_2v_{xx}(0,t)=\alpha_3w_{xx}(0,t),\; t\in (0,T),
\end{equation}
where $\alpha_2,\ \alpha_3,\ \beta_2$ and $\beta_3$ are real constants satisfying  $3\beta_i^2\leq 1$ for $i=2,3$ imply the uniqueness argument.
  
We now define the following two classes of boundary conditions that is coherent with the approach used here, which involves the particular boundary conditions \eqref{v1}-\eqref{v3} and \eqref{2v1}-\eqref{2v3}.
\begin{definition}
	Given $a_2,\ a_3,\ b_2,\ b_3,\ c_2$ and $c_3$ real constants, we call type 1  boundary conditions  for the Cauchy problem \eqref{KDV}-\eqref{initial} if these satisfy the following boundary conditions
	at the vertex:
	\begin{equation}\label{v111}
	u(0,t)=a_2v(0,t)=a_3w(0,t),\; t\in (0,T),
	\end{equation}
	\begin{equation}\label{v211}
	u_x(0,t)=b_2v_x(0,t)+b_3w_x(0,t),\; t\in (0,T)
	\end{equation}
	and
	\begin{equation}\label{v311}
	u_{xx}(0,t)=c_2v_{xx}(0,t)+c_3w_{xx}(0,t),\; t\in (0,T).
	\end{equation}
	
	\end{definition}
	
	\begin{definition}
		Given $a_2,\ a_3,\ b_2,\ b_3,\ c_2$ and $c_3$ real constants, we call type 2 boundary conditions for the Cauchy problem \eqref{KDV}-\eqref{initial} if these satisfy 
		\begin{equation}\label{v1111}
		u(0,t)=a_2v(0,t)+a_3w(0,t),\; t\in (0,T),
		\end{equation}
		\begin{equation}\label{v2111}
		u_x(0,t)=b_2v_x(0,t)+b_3w_x(0,t),\; t\in (0,T)
		\end{equation}
		and
		\begin{equation}\label{v3111}
		u_{xx}(0,t)=c_2v_{xx}(0,t)=c_3w_{xx}(0,t),\; t\in (0,T).
		\end{equation}
		
	\end{definition}


It is well-known that the trace operator $\gamma_0:u(x)\mapsto u(0)$ is well-defined on $H^{s}(\R^+)$ for $s>\frac{1}{2}$.  Hence, on the case $s>\frac12$ we will assume the following additional condition 
\begin{equation}\label{v11}
u_0(0)=a_2v_0(0)=a_3w_0(0)
\end{equation}
 for initial data for the Cauchy problem \eqref{KDV}-\eqref{initial} with type 1 boundary conditions and 
\begin{equation}\label{v12}
	u_0(0)=a_2v_0(0)+a_3w_0(0)
\end{equation}
for type 2 boundary conditions.
\begin{figure}[htp]\label{figure1}
	\centering 
	\begin{tikzpicture}[scale=3]
	\draw [thick, dashed] (-1.3,0)--(-1,0);
	\draw[thick](-1,0)--(0,0);
	\node at (-0.5,0.1){$(-\infty,0)$};
	\draw[thick](0,0)--(0.89,0.45);
	\node at (0.46,0.33)[rotate=30]{$(0,+\infty)$};
		\draw [thick, dashed] (0.89,0.45)--(1.19,0.6);
		\draw[thick](0,0)--(0.89,-0.45);
			\node at (0.46,-0.33)[rotate=-30]{$(0,+\infty)$};
			\draw [thick, dashed] (0.89,-0.45)--(1.19,-0.6);
	 \fill (0,0)  circle[radius=1pt];
	\end{tikzpicture}
	\caption{A star graph with three edges ($\mathcal{Y}$-junction)}
	\end{figure}
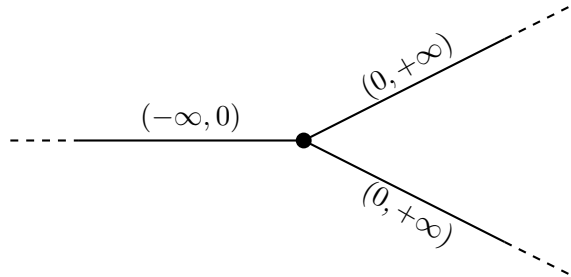

\subsection{Principal results}
Now we state our main result.

\begin{theorem}\label{theorem0}
	Let  $-\frac12<s<\frac32$, with $s\neq \frac{1}{2}$. Assume that $u_0,$ $v_0$ and $w_0$ satisfy \eqref{spacetrace}.
	\begin{itemize}
		\item [(i)]For a fixed  $s$ suppose that there exists a real constant $\lambda_i(s)$  satisfying
		\begin{equation}\label{sclarsconditions}
		 \max\{s-1,0\}<\lambda_i(s)<\min\left\{s+\frac12,\frac12\right\}\  \text{for}\ i=1,2,3,4,
		 \end{equation}
		  such that the following matrix
		\begin{equation}\label{matrix}
		\begin{split}
	\mathbf{M}(\lambda_1,\lambda_2, \lambda_3, & \lambda_4)\\
		&:=\left[\begin{array}{cccc}
		2 \text{sin}\left(\frac{\pi}{3}\lambda _1+\frac{\pi}{6}\right)&
		-a_2e^{i\pi\lambda_3}& 0&2 \text sin\left(\frac{\pi}{3}\lambda _2+\frac{\pi}{6}\right)\\
		2 \text {sin}\left(\frac{\pi}{3}\lambda _1+\frac{\pi}{6}\right)& 0&-a_3e^{i\pi\lambda_4}&2\text sin\left(\frac{\pi}{3}\lambda _2+\frac{\pi}{6}\right)\\
		2 \text {sin}\left(\frac{\pi}{3}\lambda _1-\frac{\pi}{6}\right)&-b_2e^{i\pi(\lambda_3-1)}&-b_3e^{i\pi(\lambda_4-1)}&2 \text {sin}\left(\frac{\pi}{3}\lambda _2-\frac{\pi}{6}\right)\\
		2 \text {sin}\left(\frac{\pi}{3}\lambda _1-\frac{\pi}{2}\right)&-c_2e^{i\pi(\lambda_3-2)} &-c_3e^{i\pi(\lambda_4-2)} &2 \text {sin}\left(\frac{\pi}{3}\lambda _2-\frac{\pi}{2}\right)
		\end{array}\right]
		\end{split}
		\end{equation}
		is invertible. Then there exists a positive time $T>0$ and a distributional solution $(u,v,w)$ in the space $C([0,T]:H^s(\mathcal{Y}))$, for the Cauchy problem \eqref{KDV}-\eqref{initial} with type 1 boundary conditions, satisfying the additional compatibility condition \eqref{v11} on the case $\frac12<s<\frac32$. Furthermore the data-to-solution map $(u_0,v_0,w_0)\mapsto (u,v,w)$  is locally Lipschitz continuous from $H^s(\mathcal{Y})$ to $C([0,T];H^s(\mathcal{Y}))$.
	 
		\item [(ii)] For a fixed  $s$ suppose that there exists  a real constant $\lambda_i(s)$  with $$\max\{s-1,0\}<\lambda_i(s)<\min\left\{s+\frac12,\frac12\right\}\, \text{for}\; i=1,2,3,4,$$  such that the following matrix
		\begin{equation}\label{matrix2}
		\begin{split}
		\mathbf{M}(\lambda_1,\lambda_2, \lambda_3,& \lambda_4)\\
		&:=\left[\begin{array}{cccc}
		2 \text sin\left(\frac{\pi}{3}\lambda _1+\frac{\pi}{6}\right)&
		-a_2e^{i\pi\lambda_3}& -a_3e^{i\pi\lambda_4}&2 \text sin\left(\frac{\pi}{3}\lambda _2+\frac{\pi}{6}\right)\\
		2 \text sin\left(\frac{\pi}{3}\lambda _1-\frac{\pi}{6}\right)&-b_2e^{i\pi(\lambda_3-1)}&-b_3e^{i\pi(\lambda_4-1)}&2 \text sin\left(\frac{\pi}{3}\lambda _2-\frac{\pi}{6}\right)\\
		2 \text sin\left(\frac{\pi}{3}\lambda _1-\frac{\pi}{2}\right)&-c_2e^{i\pi(\lambda_3-2)} &0&2 \text sin\left(\frac{\pi}{3}\lambda _2-\frac{\pi}{2}\right)\\
		2 \text{sin}\left(\frac{\pi}{3}\lambda _1-\frac{\pi}{2}\right)&0 &-c_3e^{i\pi(\lambda_4-2)} &2 \text sin\left(\frac{\pi}{3}\lambda _2-\frac{\pi}{2}\right)
		\end{array}\right]
		\end{split}
		\end{equation}
		is invertible. Then there exists a positive time $T>0$ and a distributional solution $(u,v,w)$ in the space $C([0,T]:H^s(\mathcal{Y}))$, for the Cauchy problem \eqref{KDV}-\eqref{initial} with type 2 boundary conditions, satisfying the additional compatibility condition \eqref{v12} on the case $\frac12<s<\frac32$. Furthermore the data-to-solution map $(u_0,v_0,w_0)\mapsto (u,v,w)$  is locally Lipschitz continuous from $H^s(\mathcal{Y})$ to $C([0,T];H^s(\mathcal{Y}))$.
		\end{itemize}	
		\end{theorem}
	
	\begin{remark}
		In Theorem \ref{theorem0} 	the indexes $\lambda_i$, $i=1,2,3,4$, are associated to the Duhamel boundary operator classes associated to linear version of KdV equation (see Subsection \ref{section4}) and depend of the regularity index $s$.
	\end{remark}
	
	\begin{remark}
		The inversion of the matrix \eqref{matrix} and \eqref{matrix2} condition is necessary in order to reformulate the Cauchy problem \eqref{KDV}-\eqref{initial} in an integral version  by using the Duhamel boundary forcing operators.
	\end{remark}

As the consequence of Theorem \ref{theorem0}, we can obtain the following result for the special boundary conditions \eqref{v1}-\eqref{v3} and \eqref{2v1}-\eqref{2v3}, which is appropriate for our formal uniqueness calculations associated to the linear version of the KdV equation on a star graph.
\begin{corollary}\label{theorem}Let  $-\frac12<s<\frac32$  with $s\neq \frac{1}{2}$ and $\alpha_2,\, \alpha_3,\, \beta_2,\, \beta_3\in \R$ satisfy $\frac{1}{\alpha_2^2}+\frac{1}{\alpha_3^2}+\frac{\beta_3}{\alpha_3}+\frac{\beta_2}{\alpha_2}\neq-1$. Assume that $u_0,$ $v_0$ and $w_0$ satisfy \eqref{spacetrace}. Then there exists a positive time $T>0$ and a distributional solution $(u,v,w)\in C([0,T];H^s(\mathcal{Y}))$ for the Cauchy problem \eqref{KDV}-\eqref{initial} with boundary conditions \eqref{v1}-\eqref{v3}, and the initial conditions satisfying additional conditions \eqref{v11} for $\frac{1}{2}<s<\frac{3}{2}$. Furthermore the data-to-solution map $(u_0,v_0,w_0)\mapsto (u,v,w)$  is locally Lipschitz continuous from $H^s(\mathcal{Y})$ to $C([0,T];H^s(\mathcal{Y}))$.
\end{corollary}
\begin{corollary}\label{theorem1} The same result of the Corollary \ref{theorem} is valid for the Cauchy problem \eqref{KDV}-\eqref{initial} with boundary conditions \eqref{2v1}-\eqref{2v3}, and the initial conditions satisfying additional condition \eqref{v12} for $\frac{1}{2}<s<\frac{3}{2}$. 
\end{corollary}

\begin{remark}Note that, in Corollary \ref{theorem} and Corollary \ref{theorem1} we don't need of the assumptions $\beta_i^2<1$, $i=2,3$, obtained in the previous formal uniqueness calculations for the associated linear problem \eqref{KDVlinear}.
	\end{remark}


The approach used to prove the main result is based on the arguments developed in \cite{Cavalcante}, \cite{CC}, \cite{CK}, \cite{Holmer1} and \cite{Holmer}. The main idea to prove Theorem \ref{theorem0} is the construction of an auxiliary forced Cauchy problem in all $\R$, analogous to the \eqref{KDV}; more precisely:
\begin{equation}\label{KDVF}
\begin{cases}
u_t+u_{xxx}+u_xu=\mathcal{T}_1(x)h_1(t)+\mathcal{T}_2(x)h_2(t),& (x,t)\in\R\times (0,T),\\
v_t+v_{xxx}+v_xv=\mathcal{T}_3(x)h_3(t),& (x,t)\in\R\times (0,T),\\
w_t+w_{xxx}+w_xw=\mathcal{T}_4(x)h_4(t),& (x,t)\in\R\times (0,T),\\
u(x,0)=\widetilde{u}_0(x),\; v(x,0)=\widetilde{v}_0(x)\; w(x,0)=\widetilde{w}_0(x),& x\in \R. 
\end{cases}
\end{equation}
where $\mathcal{T}_1$ and $\mathcal{T}_2$  are distributions supported in a positive half-line $\R^+$, $\mathcal{T}_3$ and $\mathcal{T}_4$ are distributions supported in the negative half-line $\mathbb{R}^-$,  $\widetilde{u}_0,\ \widetilde{v}_0$ and $\widetilde{w}_0$ are nice extensions of $u_0$, $v_0$ and $w_0$ in $\mathbb{R}$. The  boundary forcing functions $h_1$, $h_2$, $h_3$ and $h_4$ are selected to ensure that the vertex conditions are satisfied.

The solution of forced Cauchy problem \eqref{KDVF} satisfying the vertex conditions is constructed using the classical restricted norm method of Bourgain (see \cite{Bourgain1} and \cite{KPV}) and the inversion of a Riemann-Liouville fractional integration operator.

Following \cite{Cavalcante} and \cite{Holmer} we consider the distributions $\mathcal{T}_1=\frac{x_{-}^{\lambda_1-1}}{\Gamma(\lambda_1)}$, $\mathcal{T}_2=\frac{x_{-}^{\lambda_2-1}}{\Gamma(\lambda_2)}$,  $\mathcal{T}_3=\frac{x_{+}^{\lambda_3-1}}{\Gamma(\lambda_3)}$ and $\mathcal{T}_4=\frac{x_{+}^{\lambda_4-1}}{\Gamma(\lambda_4)}$, where 
\begin{equation}
\left\langle \frac{x_+^{\lambda-1}}{\Gamma(\lambda)},\phi\right\rangle=\int_0^{+\infty}\frac{x^{\lambda-1}}{\Gamma(\lambda)}\phi(x)dx,\; \text{for}\; \text{Re}\ \lambda>0.
\end{equation}
For other values of $\lambda$ we can define $\frac{x_+^{\lambda-1}}{\Gamma(\lambda)}=\frac{d^k}{dx}\frac{x^{\lambda+k-1}}{\Gamma(\lambda+k)}$, for any integer $k$ satisfying $k+\text{Re}\lambda>0$. Finally, we define $\frac{x_{-}^{\lambda-1}}{\Gamma(\lambda)}=e^{i\pi \lambda}\frac{(-x)_{+}^{\lambda-1}}{\Gamma(\lambda)}$.

The crucial point here are the appropriate choices of the parameters $\lambda_i$ and the functions $h_i$, for $i=1,2,3,4$, that will depend on the regularity index $s$.

\begin{remark}We believe that the same approach used to prove Theorem \ref{theorem} can provide similar results for the KdV equation in other star graphs and possibly for other nonlinear dispersive equations. For example, a treatment for the nonlinear Schr\"ondiger equation on a star graphs can be done using the classes of Duhamel boundary operators developed by Holmer \cite{Holmer1} and Cavalcante \cite{Cavalcante}. 
	\end{remark}
	
\begin{remark}
 Our results is a conditional local well-posedness, that is a weaker notion discussed by Kato \cite{Kato} in which solutions are only known to be unique if they satisfy
 additional auxiliary conditions, for a reformulation of the Cauchy problem \eqref{KDV}-\eqref{initial} for 1 and 2 types boundary conditions as an integral equation posed in $\R$, such that solves the formulation  \eqref{KDV}-\eqref{initial} on a star graph.
 \end{remark}

 We note that our best regularity attained is $s=-\frac12^+$ (see Remark \ref{remark}) and well-posedness for less regular initial data exists in the literature for the KdV equation in other domains. For example; in full real line the best regularity obtained was $s=-\frac{3}{4}$ given independently by Guo \cite{Guo} and Kishimoto \cite{Kishimoto}, in context of half-lines the best regularity was $s=-\frac{3}{4}^+$ that was obtained independently by Holmer \cite{Holmer}
 and Bona, Sun and Zhang \cite{BSZ2} in 2006.

\subsection{Organization of this paper}This paper is organized as follows: in the next section, we discuss some notation, introduce function spaces and recall some needed properties of these function spaces, and review the definition and basic properties of the Riemann-Liouville fractional integral. Sections \ref{section3},  are devoted to the needed estimates for the linear group and the Duhamel boundary forcing operators classes. In Section \ref{sectionduhamel} we state the estimates for the classical Duhamel
inhomogeneous solution operator.  Section \ref{section6} is devoted to prove Theorem \ref{theorem}. Finally, Sections \ref{section7} and \ref{section8} are devoted to proving  Corollaries \ref{theorem} and \ref{theorem1}, respectively .

\section{\textbf{Preliminaries}}

Here we introduce some notations, function spaces and the operator Riemann-Liouville fractional integral.
\subsection{Notations} For $\phi=\phi(x)\in S(\mathbb{R})$,  $\displaystyle \hat{\phi}(\xi)=\int e^{-i\xi x}\phi(x)dx$ denotes the Fourier transform of $\phi$. For $u=u(x,t)\in S(\mathbb{R}^2)$, $\hat{u}=\hat{u}(\xi,\tau)=\displaystyle \int e^{-i(\xi x+\tau t)}u(x,t)dxdt$ denotes its space-time Fourier transform, $\mathcal{F}_{x}u(\xi,t)$  its space Fourier transform and $\mathcal{F}_{t}u(x,\tau)$ its time Fourier transform. 

For any real number $\xi$ we put $\langle \xi \rangle:=1+|\xi|$ and $f(\xi,\tau) \lesssim g(\xi,\tau)$ means that there is a constant $C$ such that 
$f(\xi,\tau) \leq Cg(\xi,\tau)$ for all $(\xi, \tau)\in \R^2$. 
The characteristic function of an arbitrary set $A$ is denoted by $\chi_{A}$. Throughout the paper, we fix a cutoff function $\psi \in C_0^{\infty}(\mathbb{R})$ such that
$$
\psi(t)=
\begin{cases}
1 & \text{if}\; |t|\le 1,\\
0 & \text{if}\; |t|\ge 2
\end{cases}
$$ 
and $\mathbb{R}^*=\mathbb{R}\setminus \{0\}.$

\subsection{Function Spaces}
For $s\geq 0$ we say that $\phi \in H^s(\mathbb{R}^+)$ if exists $\tilde{\phi}\in H^s(\mathbb{R})$ such that 
$\phi=\tilde{\phi}|_{\R+}$.  In this case we set $\|\phi\|_{H^s(\mathbb{R}^+)}:=\inf\limits_{\tilde{\phi}}\|\tilde{\phi}\|_{H^{s}(\mathbb{R})}$. For $s\geq 0$ define $$H_0^s(\mathbb{R}^+)=\Big\{\phi \in H^{s}(\mathbb{R}^+);\,\text{supp} (\phi) \subset[0,+\infty) \Big\}.$$ For $s<0$, define $H^s(\mathbb{R}^+)$ and $H_0^s(\mathbb{R}^+)$  as the dual space of $H_0^{-s}(\mathbb{R}^+)$ and  $H^{-s}(\mathbb{R}^+)$, respectively. We define the usual Sobolev spaces for functions defined on the junction $\mathcal{Y}$ as 
\begin{equation}
H^s(\mathcal{Y})=H^s(\R^-)\times H^s(\R^+)\times H^s(\R^+).
\end{equation}

Also define 
$$C_0^{\infty}(\mathbb{R}^+)=\Big\{\phi\in C^{\infty}(\mathbb{R});\, \text{supp}(\phi) \subset [0,+\infty)\Big\}$$
and $C_{0,c}^{\infty}(\mathbb{R}^+)$ as those members of $C_0^{\infty}(\mathbb{R}^+)$ with compact support. We recall that $C_{0,c}^{\infty}(\mathbb{R}^+)$ is dense in $H_0^s(\mathbb{R}^+)$ for all $s\in \mathbb{R}$. A definition for $H^s(\R^-)$ and $H_0^s(\R^-)$ can be given analogous to that for $H^s(\R^+)$ and $H_0^s(\R^+)$.

The following results summarize useful properties of the Sobolev spaces on the half-line. For the proofs we refer the reader \cite{CK}.

\begin{lemma}\label{sobolevh0}
	For all $f\in H^s(\mathbb{R})$  with $-\frac{1}{2}<s<\frac{1}{2}$ we have
	\begin{equation*}
	\|\chi_{(0,+\infty)}f\|_{H^s(\mathbb{R})}\lesssim  \|f\|_{H^s(\mathbb{R})}.
	\end{equation*}
\end{lemma}
\begin{lemma}\label{alta}
	If $\frac{1}{2}<s<\frac{3}{2}$ the following statements are valid:
	\begin{enumerate}
		\item [(a)] $H_0^s(\R^+)=\big\{f\in H^s(\R^+);f(0)=0\big\},$\medskip
		\item [(b)] If  $f\in H^s(\R^+)$ with $f(0)=0$, then $\|\chi_{(0,+\infty)}f\|_{H_0^s(\R^+)}\lesssim \|f\|_{H^s(\R^+)}$.
	\end{enumerate}
\end{lemma}

\begin{lemma}\label{cut}
	If $f\in  H_0^s(\mathbb{R}^+)$ with $s\in \R$, we then have
		\begin{equation*}
	\|\psi f\|_{H_0^s(\mathbb{R}^+)}\lesssim \|f\|_{H_0^s(\mathbb{R}^+)}.
	\end{equation*}
	
\end{lemma}
\begin{remark}
	In Lemmas \ref{sobolevh0}, \ref{alta} and  \ref{cut} all the constants $c$ only depend on $s$ and $\psi$.
\end{remark}

 We denote by 
$X^{s,b}$ the so called Bourgain spaces associated to linear KdV equation; more precisely, $X^{s,b}$  is the completion of $S'(\mathbb{R}^2)$ with respect to the norm
\begin{equation*}\label{Bourgain-norm}
\|w\|_{X^{s,b}(\phi)}=\|\langle\xi\rangle^s\langle\tau-\xi^3\rangle^b\hat{w}(\xi,\tau) \|_{L_{\tau}^2L^2_{\xi}}.
\end{equation*}
To obtain our results we also need define the following auxiliary modified Bougain spaces of \cite{Holmer}. Let  $U^{s,b}$ and $V^{\alpha}$ the completion of $S'(\R^2)$ with respect to the norms:
\begin{align*}
&\|w\|_{U^{s,b}}=\left(\int\int \langle \tau\rangle^{2s/3} \langle \tau-\xi^3\rangle^{2b} |\widehat{w}(\xi,\tau)|^2d\xi d\tau\right)^{\frac{1}{2}}\label{bourgain-uxiliar-2}\\\intertext{and}
&\|w\|_{V^{\alpha}}=\left(\int\int  \langle \tau\rangle^{2\alpha} |\widehat{w}(\xi,\tau)|^2d\xi d\tau\right)^{\frac{1}{2}}.
\end{align*}

Next nonlinear estimates, in the context of the KdV equation, for $b<\frac{1}{2}$, was derived by Holmer in \cite{Holmer}.
\begin{lemma}\label{bilinear1}
	\begin{itemize}
		\item [(a)]
		Given $s>-\frac{3}{4}$, there exists $b=b(s)<\frac12$ such that for all $\alpha>\frac{1}{2}$
		we have
		\begin{equation*}
		\big\|\partial_x (v_1v_2)\big\|_{X^{s,-b}}\lesssim \|v_1\|_{X^{s,b}\cap V^{\alpha}}\|v_2\|_{X^{s,b}\cap V^{\alpha}}.
		\end{equation*}
		\item[(b)] Given $-\frac{3}{4}<s<3$, there exists $b=b(s)<\frac12$ such that for all $\alpha>\frac{1}{2}$
		we have
		\begin{equation*}
		\big\|\partial_x (v_1v_2)\big\|_{X^{s,-b}}\lesssim\|v_1\|_{X^{s,b}\cap V^{\alpha}}\|v_2\|_{X^{s,b}\cap V^{\alpha}}.
		\end{equation*}
		
	\end{itemize}
\end{lemma}

\subsection{Riemann-Liouville fractional integral}
For Re $\alpha>0$, the tempered distribution $\frac{t_+^{\alpha-1}}{\Gamma(\alpha)}$ is defined as a locally integrable function by 
\begin{equation*}
\left \langle \frac{t_+^{\alpha-1}}{\Gamma(\alpha)},\ f \right \rangle:=\frac{1}{\Gamma(\alpha)}\int_0^{+\infty} t^{\alpha-1}f(t)dt.
\end{equation*}

For Re $\alpha>0$, integration by parts implies that
\begin{equation*}
\frac{t_+^{\alpha-1}}{\Gamma(\alpha)}=\partial_t^k\left( \frac{t_+^{\alpha+k-1}}{\Gamma(\alpha+k)}\right)
\end{equation*}
for all $k\in\mathbb{N}$. This expression allows to extend the definition, in the sense of distributions,  of $\frac{t_+^{\alpha-1}}{\Gamma(\alpha)}$ to all $\alpha \in \mathbb{C}$.

If $f\in C_0^{\infty}(\mathbb{R}^+)$, we define
\begin{equation*}
\mathcal{I}_{\alpha}f=\frac{t_+^{\alpha-1}}{\Gamma(\alpha)}*f.
\end{equation*}
Thus, for Re $\alpha>0$,
\begin{equation*}
\mathcal{I}_{\alpha}f(t)=\frac{1}{\Gamma(\alpha)}\int_0^t(t-s)^{\alpha-1}f(s)ds
\end{equation*} 
and notice that 
$$\mathcal{I}_0f=f,\quad  \mathcal{I}_1f(t)=\int_0^tf(s)ds,\quad \mathcal{I}_{-1}f=f'\quad  \text{and}\quad  \mathcal{I}_{\alpha}\mathcal{I}_{\beta}=\mathcal{I}_{\alpha+\beta}.$$

The following results state important properties of the Riemann-Liouville fractional integral operator. The proof of them can be found in \cite{Holmer}.
\begin{lemma}
	If $f\in C_0^{\infty}(\mathbb{R}^+)$, then $\mathcal{I}_{\alpha}f\in C_0^{\infty}(\mathbb{R}^+)$ 
	for all $\alpha \in \mathbb{C}$.
\end{lemma}

\begin{lemma}\label{lio-lemaint}
	If $0\leq \alpha <\infty$,\, $s\in \mathbb{R}$ and $\varphi \in C_0^{\infty}(\mathbb{R})$, then we have
	\begin{align}
	&\|\mathcal{I}_{-\alpha}h\|_{H_0^s(\mathbb{R}^+)}\leq c \|h\|_{H_0^{s+\alpha}(\mathbb{R}^+)}\label{lio}\\
	\intertext{and}
	&\|\varphi \mathcal{I}_{\alpha}h\|_{H_0^s(\mathbb{R}^+)}\leq c_{\varphi} \|h\|_{H_0^{s-\alpha}(\mathbb{R}^+)}.\label{lemaint}
	\end{align}
\end{lemma}

\section{The linear versions}\label{section3}
\subsection{Linear group associated to the KdV equation}
The linear unitary  group $e^{-t\partial_x^3}:S'(\mathbb{R})\rightarrow S'(\mathbb{R})$ associated to the linear KdV equation is defined by
\begin{equation*}
e^{-t\partial_x^3}\phi(x)=\Big(e^{it\xi^3}\widehat{\phi}(\xi)\Big)^{\lor{}}(x),
\end{equation*}
that satisfies
\begin{equation}\label{lineark}
\begin{cases}
(\partial_t+\partial_x^3)e^{-t\partial_x^3}\phi(x,t)=0& \text{for}\quad (x,t)\in \mathbb{R}\times\mathbb{R},\\
e^{-t\partial_x^3}(x,0)=\phi(x)&\text{for}\quad  x\in\mathbb{R}.
\end{cases} 
\end{equation}

The next estimates were proven in \cite{Holmer}.

\begin{lemma}\label{grupok}
	Let $s\in\mathbb{R}$ and  $0< b<1$. If $\phi\in H^s(\mathbb{R})$, then we have \medskip
	\begin{enumerate}
		\item[(a)] (\textbf{space traces}) $\|e^{-t\partial_x^3}\phi(x)\|_{C\big(\mathbb{R}_t;\,H^s(\mathbb{R}_x)\big)}\lesssim \|\phi\|_{H^s(\mathbb{R})}$; \medskip
		\item[(b)] (\textbf{time traces})$\|\psi(t)\partial_x^j e^{-t\partial_x^3}\phi(x)\|_{C\big(\mathbb{R}_x;\,H^{(s+1-j)/3}(\mathbb{R}_t)\big)}\lesssim \|\phi\|_{H^s(\mathbb{R})}$ , for $j\in \mathbb{N}$; \medskip 
		\item [(c)] (\textbf{Bourgain spaces}) $\|\psi(t)e^{-t\partial_x^3}\phi(x)\|_{X^{s,b}\cap V^{\alpha}}\lesssim \|\phi\|_{H^s(\mathbb{R})}$ .
	\end{enumerate}
\end{lemma}

\begin{remark}
	The spaces $V^{\alpha}$ introduced in \cite{Holmer}  give us  useful auxiliary norms of the classical Bourgain spaces in order to validate the bilinear estimates associated to the KdV equation for $b<\frac12$ (see Lemma \ref{bilinear1}).  
\end{remark}

\subsection{The Duhamel boundary forcing operator associated to the linear KdV equation}\label{section4}
Now we give the properties of the Duhamel boundary forcing operator introduced in \cite{Holmer}, that is
\begin{equation}\label{lk}
\begin{split}
\mathcal{V}g(x,t)&=3\int_0^te^{-(t-t')\partial_x^3}\delta_0(x)\mathcal{I}_{-2/3}g(t')dt'\\
&=3\int_0^t A\left(\frac{x}{(t-t')^{1/3}}\right)\frac{\mathcal{I}_{-2/3}g(t')}{(t-t')^{1/3}}dt',
\end{split}
\end{equation}
defined for all $g\in C_0^{\infty}(\mathbb{R}^+)$ and $A$ denotes the Airy function 
$$A(x)=\frac{1}{2\pi}\int_{\xi}e^{ix\xi}e^{i\xi^3}d\xi.$$

From definition of $\mathcal{V}$  it follows that 
\begin{equation}\label{forcingk}
\begin{cases}
(\partial_t+\partial_x^3)\mathcal{V}g(x,t)=3\delta_0(x)\mathcal{I}_{-\frac{2}{3}}g(t)& \text{for}\quad (x,t)\in \mathbb{R}\times\mathbb{R},\\
\mathcal{V}g(x,0)=0& \text{for}\quad x\in\mathbb{R}.
\end{cases}
\end{equation}

The proof of the results exhibited in this section was shown in \cite{Holmer}.

\begin{lemma}\label{lemacrb1}
	Let $g\in C_0^{\infty}(\mathbb{R}^+)$ and consider a fixed time  $t\in[0,1]$. Then, 
	\begin{enumerate}
		\item[(a)]the functions $\mathcal{V}g(\cdot,t)$ and $\partial_x\mathcal{V}g(\cdot,t)$ are continuous in $x$ for all $x\in\mathbb{R}$. Moreover, they satisfy the spatial decay bounds
		$$
		|\mathcal{V}g(x,t)|+|\partial_x\mathcal{V}g(x,t)|\leq c_k\|g\|_{H^{k+1}}\langle x\rangle^{-k}\quad \text{for all}\;  k \geq 0; 
		$$
		\item[(b)] the function $\partial_x^2\mathcal{V}g(x,t)$ is continuous in $x$ for all $x\neq 0$ and has a step discontinuity
		of size $3\mathcal{I}_{\frac{2}{3}}g(t)$ at $x=0$. Also, $\partial_x^2\mathcal{V}g(x,t)$ satisfies the spatial decay bounds
		$$
		|\partial_x^2\mathcal{V}g(x,t)|\leq c_k \|f\|_{H^{k+2}}\langle x\rangle^{-k} \quad \text{for all}\; k\geq 0.
		$$
	\end{enumerate}
\end{lemma}

Since $A(0)=\frac{1}{3\Gamma\big(\frac{2}{3}\big)}$ from  \eqref{lk} we have that $\mathcal{V}g(0,t)=g(t).$

For the convenience to the reader, we present here an application of the operator $\mathcal{V}$ to solve a linear version of the IBVP associated to the KdV equation on the positive half-line, given in \cite{CK}.  Set 
\begin{equation}
v(x,t)=e^{-t\partial_x^3}\phi(x)+\mathcal{V}\big(g-e^{-\cdot\partial_x^3}\phi\big|_{x=0}\big)(x,t),
\end{equation}
where $g\in C_0^{\infty}(\R^+)$ and $\phi \in S(\R)$.

Then from \eqref{lineark} and \eqref{forcingk} we see that $v$ solves the linear problem 
\begin{equation}
\begin{cases}
(\partial_t+\partial_x^3)v(x,t)=0& \text{for}\quad (x,t)\in \R^*\times\mathbb{R},\\
v(x,0)=\phi(x)& \text{for}\quad x\in\mathbb{R},\\
v(0,t)=g(t)& \text{for}\quad t\in(0, +\infty),
\end{cases}
\end{equation}
in the sense of distributions, and then this would suffice to solve the IBVP on the right half-line associated to linear KdV equation.

Now, we consider the second boundary forcing operator associated to the linear KdV equation:
\begin{equation}
\mathcal{V}^{-1}g(x,t)=\partial_x\mathcal{V}\mathcal{I}_{\frac13}g(x,t)=3\int_0^tA'\left( \frac{x}{(t-t')^{1/3}}\right)\frac{\mathcal{I}_{-\frac13}g(t')}{(t-t')^{2/3}}dt'.
\end{equation}
From Lemma \ref{lemacrb1}, for all $g\in C_0^{\infty}(\mathbb{R}^+)$ the function $\mathcal{V}^{-1}g(x,t)$ is continuous in $x$ on $x\in\mathbb{R}$; moreover using that $A'(0)=-\frac{1}{3\Gamma(\frac{1}{3})}$ we get the relation  $\mathcal{V}^{-1}g(0,t)=-g(t).$

Also, the definition of $\mathcal{V}^{-1}g(x,t)$ allows us to ensure that
\begin{equation}
\begin{cases}
(\partial_t+\partial_x^3)\mathcal{V}^{-1}g(x,t)=3\delta_0'(x)\mathcal{I}_{-\frac{1}{3}}g(t)& \text{for}\quad (x,t)\in \mathbb{R}\times\mathbb{R},\\
\mathcal{V}^{-1}g(x,0)=0& \text{for}\quad x\in\mathbb{R},
\end{cases}
\end{equation}
in the sense of distributions.

Furthermore, Lemma \ref{lemacrb1} implies that the function $\partial_x\mathcal{V}f(x,t)$ is continuous in $x$ for all $x\in\mathbb{R}$ and, since $A'(0)=-\frac{1}{3\Gamma(\frac{1}{3})}$,
\begin{equation}
\partial_x\mathcal{V}g(0,t)=-\mathcal{I}_{-\frac{1}{3}}g(t).
\end{equation}
Also,  $\partial_x\mathcal{V}^{-1}g(x,t)=\partial_x^2\mathcal{V}\mathcal{I}_{\frac{1}{3}}g(x,t)$ is continuous in $x$ for $x\neq 0$ and has a step discontinuity of size $3\mathcal{I}_{-\frac{1}{3}}g(t)$ at $x=0$. Indeed,
\begin{eqnarray*}
	\lim_{x\rightarrow 0^{+}}\partial_x^2\mathcal{V}g(x,t)&=&-\int_{0}^{+\infty}\partial_y^3\mathcal{V}g(y,t)dy=\int_{0}^{+\infty}\partial_t\mathcal{V}g(y,t)dy\\
	&=&3\int_{0}^{+\infty}A(y)dy\int_0^t\partial_t\mathcal{I}_{-\frac{2}{3}}g(t')dt'=\mathcal{I}_{-\frac{2}{3}}g(t),
\end{eqnarray*}
then from Lemma \ref{lemacrb1} -(b) we have
\begin{equation*}
\lim_{x\rightarrow 0^{-}}\partial_x\mathcal{V}^{-1}g(x,t)=-2\mathcal{I}_{-\frac{1}{3}}g(t)\quad \text{and}\quad  \lim_{x\rightarrow 0^{+}}\partial_x\mathcal{V}^{-1}g(x,t)=\mathcal{I}_{-\frac{1}{3}}g(t).
\end{equation*}

Now, for convenience, we give an application of the operator $\mathcal{V}^{-1}$ to solve a IBVP linear associated to the KdV equation on the negative half-line with two boundary conditions given by Holmer \cite{Holmer}. Let $h_1(t)$ and $h_2(t)$ belonging to $C_0^{\infty}(\mathbb{R}^+)$ we have the relations:
\begin{eqnarray*}
	\mathcal{V}h_1(0,t)+\mathcal{V}^{-1}h_2(0,t)&=&h_1(t)-h_2(t),\\
	\lim_{x\rightarrow 0^{-}}\mathcal{I}_{\frac{1}{3}}\partial_x(\mathcal{V}h_1(x,\cdot)+\partial_x\mathcal{V}^{-1}h_2(x,\cdot))(t)&=&-h_1(t)-2h_2(t),\\
	\lim_{x\rightarrow 0^{+}}\mathcal{I}_{\frac{1}{3}}\partial_x(\mathcal{V}h_1(x,\cdot)+\partial_x\mathcal{V}^{-1}
	h_2(x,\cdot))(t)&=&-h_1(t)+h_2(t).
\end{eqnarray*}

For given $v_0(x)$, $g(t)$ and  $h(t)$ we assigned

\[ \Bigg[\begin{array}{c}
h_1 \vspace{0.3cm}\\
h_2 \end{array} \Bigg]:=\frac{1}{3}
\Bigg[\begin{array}{cc}
2  & -1 \vspace{0.3cm}\\
-1 & -1 \end{array} \Bigg]
\Bigg[\begin{array}{c}
g-e^{-\cdot\partial_x^3}v_0|_{x=0}\vspace{0.3cm}\\
\mathcal{I}_{\frac{1}{3}}\big(h-\partial_xe^{-\cdot\partial_x^3}v_0|_{x=0}\big)\end{array} \Bigg].\]
Then, taking  $v(x,t)=e^{-t\partial_x^3}v_0(x)+\mathcal{V}h_1(x,t)+\mathcal{V}^{-1}h_2(x,t)$ we get
\begin{equation}
\begin{cases}
(\partial_t+\partial_x^3)v(x,t)=0& \text{for}\quad (x,t)\in \mathbb{R}^*\times\mathbb{R},\\
v(x,0)=v_0(x)& \text{for}\quad x\in\mathbb{R},\\
v(0,t)=g(t)& \text{for}\quad t\in \R ,\\
\lim\limits_{x \rightarrow 0^{-}}\partial_xv(x,t)=h(t)& \text{for}\quad t\in \R,
\end{cases}
\end{equation}
in the sense of distributions.

\subsection{The Duhamel Boundary Forcing Operator Classes associated to linear KdV equation}

In order, to get our results in low regularity (see Remark \ref{remark}), we need to work with two classes of boundary forcing operators in order to obtain the required estimates for the second order derivative of traces. In this way, we define the generalization of operators $\mathcal{V}$ and $\mathcal{V}^{-1}$ given by Holmer \cite{Holmer}.

Let $\lambda\in \mathbb{C}$ with 
$\text{Re}\,\lambda>-3$ and $g\in C_0^{\infty}(\mathbb{R}^+)$. Define the operators
\begin{equation*}
\mathcal{V}_{-}^{\lambda}g(x,t)=\left[\frac{x_+^{\lambda-1}}{\Gamma(\lambda)}*\mathcal{V}\big(\mathcal{I}_{-\frac{\lambda}{3}}g\big)(\cdot,t)   \right](x)
\end{equation*}
and
\begin{equation*}
\mathcal{V}_{+}^{\lambda}g(x,t)=\left[\frac{x_-^{\lambda-1}}{\Gamma(\lambda)}*\mathcal{V}\big(\mathcal{I}_{-\frac{\lambda}{3}}g\big)(\cdot,t)   \right](x),
\end{equation*}
with $\frac{x_{-}^{\lambda-1}}{\Gamma(\lambda)}=e^{i\pi \lambda}\frac{(-x)_{+}^{\lambda-1}}{\Gamma(\lambda)}$. Then, 
using \eqref{forcingk} we have that
\begin{equation*}
(\partial_t+\partial_x^3)\mathcal{V}_{-}^{\lambda}g(x,t)=3\frac{x_{+}^{\lambda-1}}{\Gamma(\lambda)}\mathcal{I}_{-\frac{2}{3}-\frac{\lambda}{3}}g(t)
\end{equation*}
and
\begin{equation*}
(\partial_t+\partial_x^3)\mathcal{V}_{+}^{\lambda}g(x,t)=3\frac{x_{-}^{\lambda-1}}{\Gamma(\lambda)}\mathcal{I}_{-\frac{2}{3}-\frac{\lambda}{3}}g(t).
\end{equation*}

The following lemmas state properties of the operators classes $\mathcal{V}_{\pm}^{\lambda}$. For the proofs
we refer the reader \cite{Holmer}.
\begin{lemma}[\textbf{Spatial continuity and decay properties for} $\boldsymbol{\mathcal{V}_{\pm}^{\lambda}g(x,t)}$]\label{holmer1}
	Let $g\in C_0^{\infty}(\mathbb{R}^+)$. Then, we have
	\begin{equation*}
	\mathcal{V}_{\pm}^{\lambda-2}g=\partial_x^2\mathcal{V}^{\lambda}\mathcal{I}_{\frac{2}{3}}g,\quad \mathcal{V}_{\pm}^{\lambda-1}g=\partial_x\mathcal{V}^{\lambda}\mathcal{I}_{\frac{1}{3}}g\quad \text{and}\quad \mathcal{V}_{\pm}^{0}g=\mathcal{V}g.
	\end{equation*}
	Also, $\mathcal{V}_{\pm}^{-2}g(x,t)$ has a step discontinuity of size $3g(t)$ at $x=0$, otherwise for $x\neq 0$, $\mathcal{V}_{ \pm}^{-2}g(x,t)$ is continuous in $x$. For $\lambda>-2$, $\mathcal{V}_{\pm}^{\lambda}g(x,t)$ is continuous in $x$ for all $x\in\mathbb{R}$. For $-2\leq\lambda\leq 1$ and  $0\leq t\leq 1$, $\mathcal{V}_{- }^{\lambda}g(x,t)$ satisfies the following decay bounds:
	\begin{align*}
	&|\mathcal{V}_{-}^{\lambda}g(x,t)|\leq c_{m,\lambda,g}\langle x\rangle^{-m},\; \text{for all}\; x\leq 0\; \text{and} \;m\geq0,\\ 
	&|\mathcal{V}_{-}^{\lambda}g(x,t)|\leq c_{\lambda,g}\langle x\rangle^{\lambda-1}\; \text{for all}\;  x\geq 0.\\
	&|\mathcal{V}_{+}^{\lambda}g(x,t)|\leq c_{m,\lambda,g}\langle x\rangle^{-m},\; \text{for all}\; x\geq 0\; \text{and} \;m\geq0,\\ \intertext{and}
	&|\mathcal{V}_{+}^{\lambda}g(x,t)|\leq c_{\lambda,g}\langle x\rangle^{\lambda-1}\; \text{for all}\;  x\leq 0.
	\end{align*}
\end{lemma}

\begin{lemma}[\textbf{Values of} $\boldsymbol{\mathcal{V}_{\pm}^{\lambda}g(x,t)}$ \textbf{at} $\boldsymbol{x=0}$]\label{holmer2}
	For $\emph{Re}\,\lambda>-2$ and $g\in C_0^{\infty}(\R^+)$ we have 
	\begin{align*}
	&\mathcal{V}_{-}^{\lambda}g(0,t)=2\sin \left(\frac{\pi}{3}\lambda+\frac{\pi}{6}\right)g(t)\\ 
	\intertext{and}
	&\mathcal{V}_{+}^{\lambda}g(0,t)=e^{i\pi\lambda}g(t).
	\end{align*}
\end{lemma}
\begin{lemma}\label{cof}
	Let $s\in\mathbb{R}$.  The following estimates are ensured:\medskip
	\begin{enumerate}
		\item[(a)] (\textbf{space traces}) $\|\mathcal{V}_{\pm}^{\lambda}g(x,t)\|_{C\big(\mathbb{R}_t;\,H^s(\R_x)\big)}\lesssim \|g\|_{H_0^{(s+1)/3}(\mathbb{R}^+)}$ for all $s-\frac{5}{2}<\lambda<s+\frac{1}{2}$, $\lambda<\frac{1}{2}$  and $\emph{supp}(g)\subset[0,1]$.\medskip

		\item[(b)] (\textbf{time traces})
		$\|\psi(t)\partial_x^j\mathcal{V}_{\pm}^{\lambda}g(x,t)\|_{C\big(\mathbb{R}_x;\,H_0^{(s+1)/3}(\mathbb{R}_t^+)\big)}\lesssim c \|g\|_{H_0^{(s+1)/3}(\mathbb{R}^+)}$ for all $-2+j<\lambda<1+j$, for $j\in\{0,1,2\}$.\medskip

		\item[(c)] (\textbf{Bourgain spaces}) $
		\big\|\psi(t)\mathcal{V}_{\pm}^{\lambda}g(x,t)\big\|_{X^{s,b}\cap V^{\alpha}}\lesssim c \|g\|_{H_0^{(s+1)/3}(\mathbb{R}^+)}
		$ for all $s-1\leq \lambda<s+\frac{1}{2}$, $\lambda<\frac{1}{2}$, $\alpha\leq\frac{s-\lambda+2}{3}$ and  $0\leq b<\frac{1}{2}$.
	\end{enumerate}
\end{lemma}

\begin{remark}\label{remark} Note that for $\lambda=0$ the second derivative time traces estimate is not obtained, for this reason we need to work with the family $\mathcal{V}_{\pm}^{\lambda}$. Also note that the set of regularity where the spaces traces and Bourgain spaces estimates are valid depends of the index $\lambda$, for example, for $\lambda=0$ we have the Bourgain spaces estimates on the set $-1/2<s<1$.
\end{remark}
\section{\textbf{The Duhamel Inhomogeneous Solution Operator}}\label{sectionduhamel}

The classical inhomogeneous solution operator $\mathcal{K}$ associated to the KdV equation is given by
\begin{equation*}
\mathcal{K}w(x,t)=\int_0^te^{-(t-t')\partial_x^3}w(x,t')dt',
\end{equation*}
that satisfies
\begin{equation}\label{DK}
\begin{cases}
(\partial_t+\partial_x^3)\mathcal{K}w(x,t) =w(x,t)&\text{for}\quad  (x,t)\in\mathbb{R}\times\mathbb{R},\\
\mathcal{K}w(x,t) =0 & \text{for}\quad x\in\mathbb{R}.
\end{cases}
\end{equation}

Now, we summarize some useful estimates for the Duhamel inhomogeneous solution operators  $\mathcal{K}$  that will be used later in the proof of the main results and its proof can be seen in \cite{Holmer}.

\begin{lemma}\label{duhamelk} For all $s\in\mathbb{R}$ we have the following estimates:\medskip
	\begin{enumerate}
		\item[(a)](\textbf{space traces}) Let $-\frac{1}{2}<d<0$, then
		\begin{equation*}
		\|\psi(t)\mathcal{K}w(x,t)\|_{C\big(\mathbb{R}_t;\,H^s(\mathbb{R}_x)\big)}\lesssim\|w\|_{X^{s,d}}.
		\end{equation*}\medskip
		
		\item[(b)](\textbf{time traces})
		Let $-\frac{1}{2}<d<0$ and $j\in \{0,1,2\}$, then
		\begin{equation*}
		\left\|\psi(t)\partial_x^j\mathcal{K}w(x,t)\right\|_{C\big(\mathbb{R}_x;\,H^{(s+1)/3}(\mathbb{R}_t)\big)}\lesssim \begin{cases}
		\|w\|_{X^{s,d}} & \text{if}\; -1+j\leq s \leq \frac{1}{2}+j,\\
		\|w\|_{X^{s,d}}+\|w\|_{U^{s,d}}& \text{for all}\;  s\in\mathbb{R}.
		\end{cases}
		\end{equation*}\medskip

		\item[(c)](\textbf{Bourgain spaces estimates})
		Let $0<b<\frac{1}{2}$ and $\alpha>1-b$, then
		\begin{equation*}
		\|\psi(t)\mathcal{K}w(x,t)\|_{X^{s,b}\cap V^{\alpha}}\lesssim \|w\|_{X^{s,-b}}.
		\end{equation*}
	\end{enumerate}
\end{lemma}

\begin{remark}
	We note that the  time-adapted Bourgain spaces $U^{k,d}$ used in Lemma \ref{duhamelk}-(c)-(d) are introduced in order to cover the full values of regularity $s$. 
\end{remark}

\section{Proof of Theorem \ref{theorem0}}\label{section6}

Here we show the proof of the main result announced of this work. We only prove the part (i) of Theorem \ref{theorem0}, since the proof of part (ii) is very similar. We follow closely the arguments in \cite{Holmer} (see also \cite{Cavalcante}  and \cite{CC}). The proof will be divided into five steps.  

\textbf{Step 1. We will first obtain an integral equation that solves Cauchy problem \eqref{KDV}-\eqref{initial}, with type 1 boundary conditions, satisfying \eqref{v11} for $\frac12<s<\frac32$.}

We start rewriting the vertex conditions \eqref{v111}, \eqref{v211} and \eqref{v311} in terms of matrices:
\begin{equation}\label{m1}
\left[\begin{array}{ccc}
1&
-a_2& 0\\
1& 0&-a_3\\
0&0&0\\
0&0&0
\end{array}\right]\left[\begin{array}{r}
u(0,t)\\
v(0,t)\\
w(0,t)
\end{array}\right]=0,
\end{equation}
\begin{equation}\label{m2}
\left[\begin{array}{ccc}
0& 0& 0\\
0& 0&0\\
1&-b_2&-b_3\\
0&0&0
\end{array}\right]\left[\begin{array}{r}
u_x(0,t)\\
v_x(0,t)\\
w_x(0,t)
\end{array}\right]=0
\end{equation}
and
\begin{equation}\label{m3}
\left[\begin{array}{ccc}
0& 0& 0\\
0& 0&0\\
0&0&0\\
1&-c_2&-c_3
\end{array}\right]\left[\begin{array}{r}
u_{xx}(0,t)\\
v_{xx}
(0,t)\\
w_{xx}(0,t)
\end{array}\right]=0.
\end{equation}
Let $\widetilde{u}_0, \widetilde{v}_0$ and $\widetilde{w}_0$ nice extensions of $u_0,v_0$ and $w_0$, respectively satisfying
\begin{equation*}
\|\widetilde{u}_0\|_{H^s(\R)}\leq c\|u_0\|_{H^s(\R^{+})},\; \|\widetilde{v}_0\|_{H^s(\R)}\leq c\|v_0\|_{H^s(\R^{+})}\;\text{and}\; \|\widetilde{w}_0\|_{H^s(\R)}\leq c\|w_0\|_{H^s(\R^{+})}.
\end{equation*}
Initially, we look for solutions in the form
\begin{align*}
&u(x,t)=\mathcal{V}_-^{\lambda_1}\gamma_1(x,t)+\mathcal{V}_-^{\lambda_2}\gamma_2(x,t)+F_1(x,t),\\
&v(x,t)=\mathcal{V}_+^{\lambda_3}\gamma_3(x,t)+F_2(x,t),\\
&w(x,t)=\mathcal{V}_+^{\lambda_4}\gamma_4(x,t)+F_3(x,t),
\end{align*}
where $\gamma_i$ ($i=1,2,3,4$) are unknown functions and
\begin{align*}
&F_1(x,t)=e^{it\partial_x^3}\widetilde{u}_0+\mathcal{K}(uu_x)(x,t),\\ &F_2(x,t)=e^{it\partial_x^3}\widetilde{v}_0+\mathcal{K}(vv_x)(x,t),\\ &F_3(x,t)=e^{it\partial_x^3}\widetilde{w}_0+\mathcal{K}(ww_x)(x,t).
\end{align*}
By using Lemma \ref{holmer2} we see that
\begin{align}
&u(0,t)=2\ \text sin\left(\frac{\pi}{3}\lambda _1+\frac{\pi}{6}\right)\gamma_1(t)+2\  \text sin\left(\frac{\pi}{3}\lambda _2+\frac{\pi}{6}\right)\gamma_2(t)+F_1(0,t),\label{tec1}\\
&v(0,t)=e^{i\pi\lambda_3}\gamma_3(t)+F_2(0,t),\label{tec2}\\
&w(0,t)=e^{i\pi\lambda_4}\gamma_4(t)+F_3(0,t).\label{tec3}
\end{align}
Now we calculate the traces of first derivative functions. By Lemmas \ref{holmer1} and \ref{holmer2} we see that
\begin{align}
&u_x(0,t)=2\ \text sin\left(\frac{\pi}{3}\lambda _1-\frac{\pi}{6}\right)\mathcal{I}_{-\frac{1}{3}}\gamma_1(t)+2\  \text sin\left(\frac{\pi}{3}\lambda _2-\frac{\pi}{6}\right)\mathcal{I}_{-\frac{1}{3}}\gamma_2(t)+\partial_xF_1(0,t),\label{tec4}\\
&v_x(0,t)=e^{i\pi(\lambda_3-1)}\mathcal{I}_{-\frac13}\gamma_3(t)+\partial_xF_2(0,t),\label{tec5}\\
&w_x(0,t)=e^{i\pi(\lambda_4-1)}\mathcal{I}_{-\frac13}\gamma_4(t)+\partial_xF_3(0,t).\label{tec6}
\end{align}
In the same way, we calculate the traces of second derivatives functions,
\begin{align}
&u_{xx}(0,t)=2 \text sin\left(\frac{\pi}{3}\lambda _1-\frac{\pi}{2}\right)\mathcal{I}_{-\frac{2}{3}}\gamma_1(t)+2\  \text sin\left(\frac{\pi}{3}\lambda _2-\frac{\pi}{2}\right)\gamma_2(t)+\partial_x^2F_1(0,t),\label{tec7}\\
&v_{xx}(0,t)=e^{i\pi(\lambda_3-2)}\mathcal{I}_{-\frac23}\gamma_3(t)+\partial_x^2F_2(0,t),\label{tec8}\\
&w_{xx}(0,t)=e^{i\pi(\lambda_4-2)}\mathcal{I}_{-\frac23}\gamma_4(t)+\partial_x^2F_3(0,t).\label{tec9}
\end{align}
Note that by Lemmas \ref{holmer1} and \ref{holmer2} these calculus are valid for Re $\lambda>0$.

By substituting \eqref{tec1}, \eqref{tec2} and \eqref{tec3} into \eqref{m1}; \eqref{tec3}, \eqref{tec4} and \eqref{tec5} into \eqref{m2}, and \eqref{tec7}, \eqref{tec8} and \eqref{tec9} into \eqref{m3} we see that the functions $\gamma_i$ and indexes $\lambda_i$, for $i=1,2,3,4$, satisfy the expressions
\begin{equation}
\begin{split}
&\left[\begin{array}{ccc}
1&
-a_2& 0\\
1& 0&-a_3\\
0&0&0\\
0&0&0.
\end{array}\right]\left[\begin{array}{cccc}
2\ \text sin\left(\frac{\pi}{3}\lambda _1+\frac{\pi}{6}\right)&
0& 0&2\ \text sin\left(\frac{\pi}{3}\lambda _2+\frac{\pi}{6}\right)\\
0& e^{i\pi\lambda_3}&0&0\\
0&0&e^{i\pi\lambda_4}&0
\end{array}\right]\left[\begin{array}{r}
\gamma_1\\
\gamma_3\\
\gamma_4\\
\gamma_2
\end{array}\right]\\
&\;\;\;\;\;\;\;\;\;\;\;\;\;\;\;\;\;\;\;\;\;\;\;\;\;\;\;\;\;\;\;\;\;\;\;\;\;\;\;\;\;\;\;\;\;\;\;\;\;\;\;\;\;\;\;\;=-\left[\begin{array}{ccc}
1&
-a_2& 0\\
1& 0&-a_3\\
0&0&0\\
0&0&0.
\end{array}\right]\left[\begin{array}{r}
F_1(0,t)\\
F_2(0,t)\\
F_3(0,t)
\end{array}\right],
\end{split}
\end{equation}
\begin{equation}
\begin{split}
&\left[\begin{array}{ccc}
0&
0& 0\\
0& 0&0\\
1&-b_2&-b_3\\
0&0&0.
\end{array}\right]\left[\begin{array}{cccc}
2\ \text sin\left(\frac{\pi}{3}\lambda _1-\frac{\pi}{6}\right)&
0& 0&2\  \text sin\left(\frac{\pi}{3}\lambda _2-\frac{\pi}{6}\right)\\
0& e^{i(\pi\lambda_3-1)}&0&0\\
0&0&e^{i(\pi\lambda_4-1)}&0
\end{array}\right]\left[\begin{array}{r}
\gamma_1\\
\gamma_3\\
\gamma_4\\
\gamma_2
\end{array}\right]\\
&\;\;\;\;\;\;\;\;\;\;\;\;\;\;\;\;\;\;\;\;\;\;\;\;\;\;\;\;\;\;\;\;\;\;\;\;\;\;\;\;\;\;\;\;\;\;\;\;\;\;\;\;\;\;\;\;=-\left[\begin{array}{ccc}
0&
0& 0\\
0& 0&0\\
1&-b_2&-b_3\\
0&0&0.
\end{array}\right]\left[\begin{array}{r}
\partial_x\mathcal{I}_{\frac{1}{3}}F_1(0,t)\\
\partial_x\mathcal{I}_{\frac{1}{3}}F_2(0,t)\\
\partial_x\mathcal{I}_{\frac{1}{3}}F_3(0,t)
\end{array}\right]
\end{split}
\end{equation}
and
\begin{equation}
\begin{split}
&\left[\begin{array}{ccc}
0&
0& 0\\
0& 0&0\\
0&0&0\\
1&-c_2
&-c_3
\end{array}\right]\left[\begin{array}{cccc}
2\ \text sin\left(\frac{\pi}{3}\lambda _1-\frac{\pi}{2}\right)&
0& 0&2\ \text sin\left(\frac{\pi}{3}\lambda _2-\frac{\pi}{2}\right)\\
0& e^{i(\pi\lambda_3-2)}&0&0\\
0&0&e^{i(\pi\lambda_4-2)}&0
\end{array}\right]\left[\begin{array}{r}
\gamma_1\\
\gamma_3\\
\gamma_4\\
\gamma_2
\end{array}\right]\\
&\;\;\;\;\;\;\;\;\;\;\;\;\;\;\;\;\;\;\;\;\;\;\;\;\;\;\;\;\;\;\;\;\;\;\;\;\;\;\;\;\;\;\;\;\;\;\;\;\;\;\;\;\;\;\;\;=-\left[\begin{array}{ccc}
0&
0& 0\\
0& 0&0\\
0&0&0\\
1&-c_2&-c_3.
\end{array}\right]\left[\begin{array}{r}
\partial_x^2\mathcal{I}_{\frac{2}{3}}F_1(0,t)\\
\partial_x^2\mathcal{I}_{\frac{2}{3}}F_2(0,t)\\
\partial_x^2\mathcal{I}_{\frac{2}{3}}F_3(0,t)
\end{array}\right].
\end{split}
\end{equation}
It follows that,
\begin{equation}\label{m111}
\begin{split}
&\left[\begin{array}{cccc}
2\ \text sin\left(\frac{\pi}{3}\lambda _1+\frac{\pi}{6}\right)&
-a_2e^{i\pi\lambda_3}& 0&2\ \text sin\left(\frac{\pi}{3}\lambda _2+\frac{\pi}{6}\right)\\
2 \text sin\left(\frac{\pi}{3}\lambda _1+\frac{\pi}{6}\right)& 0&-a_3e^{i\pi\lambda_4}&2\text sin\left(\frac{\pi}{3}\lambda _2+\frac{\pi}{6}\right)\\
0&0&0&0\\
0&0 &0 &0
\end{array}\right]\left[\begin{array}{r}
\gamma_1\\
\gamma_3\\
\gamma_4\\
\gamma_2
\end{array}\right]\\
&=-\left[\begin{array}{c}
F_1(0,t)-a_2F_2(0,t)\\
F_1(0,t)-a_3F_3(0,t)\\
0\\
0
\end{array}\right].
\end{split}
\end{equation}

\begin{equation}\label{m222}
\begin{split}
&\left[\begin{array}{cccc}
0&
0& 0&0\\
0& 0&0&0\\
2\ \text sin\left(\frac{\pi}{3}\lambda _1-\frac{\pi}{6}\right)&-b_2e^{i\pi(\lambda_3-1)}&-b_3e^{i\pi(\lambda_4-1)}&2\  \text sin\left(\frac{\pi}{3}\lambda _2-\frac{\pi}{6}\right)\\
0&0 &0 &0
\end{array}\right]\left[\begin{array}{r}
\gamma_1\\
\gamma_3\\
\gamma_4\\
\gamma_2
\end{array}\right]\\
&=-\left[\begin{array}{c}
0\\
0\\
\partial_x\mathcal{I}_{\frac{1}{3}}F_1(0,t)-b_2\partial_x\mathcal{I}_{\frac{1}{3}}F_2(0,t)-b_3\partial_x\mathcal{I}_{\frac{1}{3}}F_3(0,t)\\
0
\end{array}\right]
\end{split}
\end{equation}
and
\begin{equation}\label{m333}
\begin{split}
&\left[\begin{array}{cccc}
0&
0& 0&0\\
0& 0&0&0\\
0&0&0&0\\
2\ \text sin\left(\frac{\pi}{3}\lambda _1-\frac{\pi}{2}\right)&-c_2e^{i\pi(\lambda_3-2)} &-c_3e^{i\pi(\lambda_4-2)} &2\  \text sin\left(\frac{\pi}{3}\lambda _2-\frac{\pi}{2}\right)
\end{array}\right]\left[\begin{array}{r}
\gamma_1\\
\gamma_3\\
\gamma_4\\
\gamma_2
\end{array}\right]\\
&=-\left[\begin{array}{c}
0\\
0\\
0\\
\partial_x^2\mathcal{I}_{\frac{2}{3}}F_1(0,t)-c_2\partial_x^2\mathcal{I}_{\frac{2}{3}}F_2(0,t)-c_3\partial_x^2\mathcal{I}_{\frac{2}{3}}F_3(0,t)
\end{array}\right].
\end{split}
\end{equation}

From \eqref{m111}, \eqref{m222} and \eqref{m333} we need to get functions $\gamma_i$ ($i=1,2,3,4$) and parameters $\lambda_i$ $(i=1,2,3,4)$ satisfying
\begin{equation}\label{M}
\begin{split}
&\left[\begin{array}{cccc}
2\ \text sin\left(\frac{\pi}{3}\lambda _1+\frac{\pi}{6}\right)&
-a_2e^{i\pi\lambda_3}& 0&2\ \text sin\left(\frac{\pi}{3}\lambda _2+\frac{\pi}{6}\right)\\
2\ \text sin\left(\frac{\pi}{3}\lambda _1+\frac{\pi}{6}\right)& 0&-a_3e^{i\pi\lambda_4}&2\ \text sin\left(\frac{\pi}{3}\lambda _2+\frac{\pi}{6}\right)\\
2\ \text sin\left(\frac{\pi}{3}\lambda _1-\frac{\pi}{6}\right)&-b_2e^{i\pi(\lambda_3-1)}&-b_3e^{i\pi(\lambda_4-1)}&2\ \text sin\left(\frac{\pi}{3}\lambda _2-\frac{\pi}{6}\right)\\
2\ \text sin\left(\frac{\pi}{3}\lambda _1-\frac{\pi}{2}\right)&-c_2e^{i\pi(\lambda_3-2)} &-c_3e^{i\pi(\lambda_4-2)} &2\ \text sin\left(\frac{\pi}{3}\lambda _2-\frac{\pi}{2}\right)
\end{array}\right]\left[\begin{array}{r}
\gamma_1\\
\gamma_3\\
\gamma_4\\
\gamma_2
\end{array}\right]\\
&=-\left[\begin{array}{c}
F_1(0,t)-a_2F_2(0,t)\\
F_1(0,t)-a_3F_3(0,t)\\
\partial_x\mathcal{I}_{\frac{1}{3}}F_1(0,t)-b_2\partial_x\mathcal{I}_{\frac{1}{3}}F_2(0,t)-b_3\partial_x\mathcal{I}_{\frac{1}{3}}F_3(0,t)\\
\partial_x^2\mathcal{I}_{\frac{2}{3}}F_1(0,t)-c_2\partial_x^2\mathcal{I}_{\frac{2}{3}}F_2(0,t)-c_3\partial_x^2\mathcal{I}_{\frac{2}{3}}F_3(0,t)
\end{array}\right].
\end{split}
\end{equation}

We denote a simplified notation of \eqref{M} as
\begin{equation}
\mathbf{M}(\lambda_1,\lambda_2, \lambda_3, \lambda_4)\boldsymbol{\gamma}=\mathbf{F},
\end{equation}
where $\mathbf{M}(\lambda_1,\lambda_2, \lambda_3, \lambda_4)$ is the first matrix that appears in \eqref{M}, $\boldsymbol{\gamma}$ is the  matrix column given by vector $(\gamma_1,\gamma_2,\gamma_3,\gamma_4)$ and $\mathbf{F}$ is the last matrix in $\eqref{M}$. By using the hypothesis of Theorem \ref{theorem0} we fix parameters $\lambda_i$, for $i=1,2,3,4$ such that
 \begin{equation}
\max\{s-1,0\}<\lambda_i(s)<\min\left\{s+\frac12,\frac12\right\}.
 \end{equation}
and the matrix $\mathbf{M}(\lambda_1,\lambda_2, \lambda_3, \lambda_4)$ is invertible.
\medskip

\textbf{Step 2. We will define the truncated integral operator and the appropriate functions space.}

Given $s$ as in the hypothesis of Theorem \ref{theorem0} we fix the parameters $\lambda_i$ and the functions $\gamma_i$ $(i=1,2,3,4)$ chosen as in the Step 1. Let $b=b(s)<\frac12$ and $\alpha(b)>1/2$ such that the estimates given in Lemma \ref{bilinear1} are valid.

Define the operator
\begin{equation}
\Lambda=(\Lambda_1,\Lambda_2,\Lambda_3)
\end{equation}
where
\begin{align*}
&\Lambda_1 u(x,t)=\psi(t)\mathcal{V}_-^{\lambda_1}\gamma_1(x,t)+\psi(t)\mathcal{V}_-^{\lambda_2}\gamma_2(x,t)+F_1(x,t),\\
&\Lambda_2v(x,t)=\psi(t)\mathcal{V}_+^{\lambda_3}\gamma_3(x,t)+F_2(x,t),\\
&\Lambda_3w(x,t)=\psi(t)\mathcal{V}_+^{\lambda_4}\gamma_4(x,t)+F_3(x,t),
\end{align*}
where
\begin{align*}
&F_1(x,t)=\psi(t)(e^{it\partial_x^3}\widetilde{u}_0+\mathcal{K}(uu_x)(x,t)),\\ &F_2(x,t)=\psi(t)(e^{it\partial_x^3}\widetilde{v}_0+\mathcal{K}(vv_x)(x,t)),\\ &F_3(x,t)=\psi(t)(e^{it\partial_x^3}\widetilde{w}_0+\mathcal{K}(ww_x)(x,t)).
\end{align*}

We consider $\Lambda$ on the Banach space $Z(s)=Z_1(s)\times Z_2(s)\times Z_3(s)$, where
\begin{equation*}
\begin{split}
Z_i(s)=\{w\in C(\R_t;&H^s(\R_x)) \cap C(\R_x;H^{\frac{s+1}{3}}(\R_t))\cap X^{s,b}\cap V^{\alpha};\\
& w_x\in C(\R_x;H^{\frac{s}{3}}(\R_t)),w_{xx}\in C(\R_x;H^{\frac{s-1}{3}}(\R_t)) \}\; (i=1,2,3),
\end{split}
\end{equation*}
with norm
\begin{equation*}
\|(u,v,w)\|_{Z(s)}=\|u\|_{Z_1(s)}+\|v\|_{Z_2(s)}+\|w\|_{Z_3(s)},
\end{equation*}
where
\begin{equation}
\begin{split}
\|u\|_{Z_i(s)}=&\|u\|_{C(\R_t;H^s(\R_x))}+\|u\|_{C(\R_x;H^{\frac{s+1}{3}}(\R_t))}+\|u\|_{ X^{s,b}}+\|u\|_{V^{\alpha}}\\
&+\|u_x\|_{C(\R_x;H^{\frac{s}{3}}(\R_t))}+\|u_{xx}\|_{C(\R_x;H^{\frac{s-1}{3}}(\R_t))}.
\end{split}
\end{equation}

\textbf{Step 3. We will prove that the functions} $\mathbf{\mathcal{V}_{-}^{\lambda_1}\gamma_1(x,t)},\ \mathbf{\mathcal{V}_{-}^{\lambda_2}\gamma_2(x,t)},\ \mathbf{\mathcal{V}_{+}^{\lambda_3}\gamma_1(x,t)}$ \textbf{and} $ \mathbf{\mathcal{V}_{-}^{\lambda_4}\gamma_4(x,t)}$ \textbf{are well defined.}

By Lemma \eqref{cof} it suffices to show that these functions are in the closure of the spaces $C_0^{\infty}(\R^+)$. By using expression \eqref{M} we see that the functions $\gamma_i$ $(i=1,2,3,4)$ are linear combinations of the functions $F_1(0,t)-a_2F_2(0,t),F_1(0,t)-a_3F_3(0,t), \partial_x\mathcal{I}_{\frac{1}{3}}F_1(0,t)-b_2\partial_x\mathcal{I}_{\frac{1}{3}}F_2(0,t)-b_3\partial_x\mathcal{I}_{\frac{1}{3}}F_3(0,t)$ and
$\partial_x^2\mathcal{I}_{\frac{2}{3}}F_1(0,t)-c_2\partial_x^2\mathcal{I}_{\frac{2}{3}}F_2(0,t)-c_3\partial_x^2\mathcal{I}_{\frac{2}{3}}F_3(0,t).$ Thus, we need to show that the functions $F_i(0,t),\  \partial_x\mathcal{I}_{\frac{1}{3}}F_i(0,t),\  \partial_x^2\mathcal{I}_{\frac{2}{3}}F_i(0,t)$ are in appropriate spaces. By using Lemmas  \ref{grupok}, \ref{cof}, \ref{duhamelk} and \ref{bilinear1} we obtain
\begin{equation}
\|F_1(0,t)\|_{H^{\frac{s+1}{3}}(\R^+)}\leq c( \|u_0\|_{H^s(\R^+)}+\|u\|_{X^{s,b}}^2+\|u\|_{Y^{\alpha}}^2),
\end{equation}
\begin{equation}
\|F_2(0,t)\|_{H^{\frac{s+1}{3}}(\R^+)}\leq c( \|v_0\|_{H^s(\R^+)}+\|v\|_{X^{s,b}}^2+\|v\|_{Y^{\alpha}}^2),
\end{equation}
\begin{equation}
\|F_3(0,t)\|_{H^{\frac{s+1}{3}}(\R^+)}\leq c( \|w_0\|_{H^s(\R^+)}+\|w\|_{X^{s,b}}^2+\|w\|_{Y^{\alpha}}^2).
\end{equation}
If $-\frac{1}{2}<s<\frac{1}{2}$ we have that $\frac16 <\frac{s+1}{3}<\frac12$. Thus Lemma \ref{sobolevh0} implies that $H^{\frac{s+1}{3}}(\R^+)=H_0^{\frac{s+1}{3}}(\R^+)$. It follows that $F_i(0,t)\in H_0^{\frac{s+1}{3}}(\R^+)$ $(\text{for}\ i=1,2,3)$ for $-\frac{1}{2}<s<\frac{1}{2}.$

If $\frac{1}{2}<s<\frac{3}{2}$, then $\frac12<\frac{s+1}{3}<\frac56$. Using the compatibility condition \eqref{v11} we have that 
\begin{equation*}
F_1(0,0)-a_2 F_2(0,0)=u(0,0)-a_2v(0,0)=u_0(0)-a_2v_0(0)=0,
\end{equation*}
\begin{equation*}
F_1(0,0)-a_3 F_3(0,0)=u(0,0)-a_3w(0,0)=u_0(0)-a_3w_0(0)=0.
\end{equation*}
Then Lemma \ref{alta} implies 
\begin{equation}\label{trace1}
\begin{split}
F_1(0,t)-a_22F_2(0,t)\in H_0^{\frac{s+1}{3}}(\R^+),\\
F_1(0,t)-a_3F_3(0,t)\in H_0^{\frac{s+1}{3}}(\R^+)
\end{split}
\end{equation}

Now using Lemmas  \ref{grupok}, \ref{cof}, \ref{duhamelk} and \ref{bilinear1} we see that
\begin{equation*}
\|\partial_x F_1(0,t)\|_{H^{\frac{s}{3}}(\R^+)}\leq c( \|u_0\|_{H^s(\R^+)}+\|u\|_{X^{s,b}}^2+\|u\|_{Y^{\alpha}}^2), 
\end{equation*}
\begin{equation*}
\|\partial_xF_2(0,t)\|_{H^{\frac{s}{3}}(\R^+)}\leq c( \|v_0\|_{H^s(\R^+)}+\|v\|_{X^{s,b}}^2+\|v\|_{Y^{\alpha}}^2), 
\end{equation*}
\begin{equation*}
\|\partial_xF_3(0,t)\|_{H^{\frac{s}{3}}(\R^+)}\leq c( \|w_0\|_{H^s(\R^+)}+\|w\|_{X^{s,b}}^2+\|w\|_{Y^{\alpha}}^2).
\end{equation*}
Since $-\frac12<s<\frac32$ we have $-\frac16<\frac{s}{3}<\frac12$, then Lemma \ref{sobolevh0} implies that the functions $\partial_xF_i(0,t)\in H_0^{\frac{s}{3}}(\R^+)$, for $i=1,2,3,4$. Then using Lemma \ref{lio} we have that
\begin{equation*}
\|\partial_x \mathcal{I}_{\frac13}F_1(0,t)\|_{H_0^{\frac{s+1}{3}}(\R^+)}\leq c( \|u_0\|_{H^s(\R^+)}+\|u\|_{X^{s,b}}^2+\|u\|_{Y^{\alpha}}^2), 
\end{equation*}
\begin{equation*}
\|\partial_x\mathcal{I}_{\frac13}F_2(0,t)\|_{H_0^{\frac{s+1}{3}}(\R^+)}\leq c( \|v_0\|_{H^s(\R^+)}+\|v\|_{X^{s,b}}^2+\|v\|_{Y^{\alpha}}^2), 
\end{equation*}
\begin{equation*}
\|\partial_x\mathcal{I}_{\frac13}F_3(0,t)\|_{H_0^{\frac{s+1}{3}}(\R^+)}\leq c( \|w_0\|_{H^s(\R^+)}+\|w\|_{X^{s,b}}^2+\|w\|_{Y^{\alpha}}^2).
\end{equation*}

Thus, we have 
\begin{equation}\label{trace2}
\partial_x\mathcal{I}_{\frac{1}{3}}F_1(0,t)-b_2\partial_x\mathcal{I}_{\frac{1}{3}}F_2(0,t)-b_3\partial_x\mathcal{I}_{\frac{1}{3}}F_3(0,t)\in H_0^{\frac{s+1}{3}}(\R^+).
\end{equation}

In the same way we can obtain
\begin{equation*}
\|\partial_x^2\mathcal{I}_{\frac{2}{3}}F_1(0,t)\|_{H_0^{\frac{s+1}{3}}(\R^+)}\leq c( \|u_0\|_{H^s(\R^+)}+\|u\|_{X^{s,b}}^2+\|u\|_{Y^{\alpha}}^2), 
\end{equation*}
\begin{equation*}
\|\partial_x^2\mathcal{I}_{\frac{2}{3}}F_2(0,t)\|_{H_0^{\frac{s+1}{3}}(\R^+)}\leq c( \|v_0\|_{H^s(\R^+)}+\|v\|_{X^{s,b}}^2+\|v\|_{Y^{\alpha}}^2), 
\end{equation*}
\begin{equation*}
\|\partial_x^2\mathcal{I}_{\frac{2}{3}}F_3(0,t)\|_{H_0^{\frac{s+1}{3}}(\R^+)}\leq c( \|w_0\|_{H^s(\R^+)}+\|w\|_{X^{s,b}}^2+\|w\|_{Y^{\alpha}}^2).
\end{equation*}
It follows that
\begin{equation}\label{trace3}
\partial_x^2\mathcal{I}_{\frac{2}{3}}F_1(0,t)-c_2\partial_x^2\mathcal{I}_{\frac{2}{3}}F_2(0,t)-c_3\partial_x\mathcal{I}_{\frac{2}{3}}F_3(0,t)\in H_0^{\frac{s-1}{3}}(\R^+).
\end{equation}
Thus, \eqref{trace1}, \eqref{trace2} and \eqref{trace3} imply that the functions $\mathcal{V}_{-}^{\lambda_1}\gamma_1(x,t),\ \mathcal{V}_{-}^{\lambda_2}\gamma_2(x,t),\ \mathcal{V}_{+}^{\lambda_3}\gamma_1(x,t)$ and $ \mathcal{V}_{-}^{\lambda_4}\gamma_4(x,t)$ are well defined.

\textbf{Step 4. We will obtain a fixed point for} $\mathbf{\Lambda}$ \textbf{in a ball of} $\mathbf{Z}$.

Using Lemmas  \ref{lio}, \ref{grupok}, \ref{cof}, \ref{duhamelk} and \ref{bilinear1} we obtain
\begin{equation}
\begin{split}
\|\Lambda(u_2,v_2,w_2)-\Lambda&(u_1,v_1,w_1)\|_{Z}\\
&\leq c(\|(u_2,v_2,w_2)\|_Z+\|(u_1,v_1,w_1)\|_{Z})\|(u_2,v_2,w_2)-(u_1,v_1,w_1)\|_{Z}
\end{split}
\end{equation}
and 
\begin{equation}
\begin{split}
\|\Lambda(u,v,w)\|_{Z}\leq c& (\|u_0\|_{H^s(\R^+)}+\|v_0\|_{H^s(\R^+)}+\|w_0\|_{H^s(\R^+)}\\
&+\|u\|_{X^{s,b}}^2+\|u\|_{Y^{\alpha}}^2+\|v\|_{X^{s,b}}^2+\|v\|_{Y^{\alpha}}^2+\|w\|_{X^{s,b}}^2+\|w\|_{Y^{\alpha}}^2).
\end{split}
\end{equation}
By taking $\|u_0\|_{H^s(\R^+)}+\|v_0\|_{H^s(\R^+)}+\|w_0\|_{H^s(\R^+}<\delta$ for $\delta>0$ suitable small, we obtain a fixed point $\Lambda(\widetilde{u},\widetilde{v},\widetilde{w})=(\widetilde{u},\widetilde{v},\widetilde{w})$ in a ball $$B=\{(u,v,w)\in Z, \|(u,v,w)\|_Z\leq 2c\delta\}.$$ It follows that the restriction
\begin{equation}
(u,v,w)=(\widetilde{u}\big|_{\R^-\times (0,1)},\widetilde{v}\big|_{\R^+\times (0,1)},\widetilde{w}\big|_{\R^+\times (0,1)})
\end{equation}
solves the Cauchy problem \eqref{KDV}-\eqref{initial} with 1 boundary conditions in the sense of distributions.

Existence of solutions for any data in $H^s(\mathcal{Y})$ follows by the standard scaling argument. Suppose we are given data $\widetilde{u}_0, \widetilde{v}_0$ and $\widetilde{w}_0$ with arbitrary size for the Cauchy problem \eqref{KDV}-\eqref{initial} with 1 boundary conditions. For $\lambda<<1$ (to be selected after) define $u_0(x)=\lambda^{2}\widetilde{u}_0(\lambda x),\ v_0(x)=\lambda^{2}\widetilde{v}_0(\lambda x)$ and $w_0(x)=\lambda^{2}\widetilde{w}_0(\lambda x)$. Taking $\lambda$ sufficiently small so that
\begin{equation}
\|u_0\|_{H^s}+\|v_0\|_{H^s}+\|w_0\|_{H^s}\leq \lambda^{\frac32}(1+\lambda^{s})(\|\widetilde{u}_0\|_{H^s(\R^+)}+\|\widetilde{v}_0\|_{H^s(\R^+)}+\|\widetilde{w}_0\|_{H^s(\R^+)})< \delta.
\end{equation}
Then using the previous argument, there exists a solution $u(x,t)$ for the problem \eqref{KDV}-\eqref{initial}, with type 1 boundary conditions, on $0\leq t\leq 1$. Then $\widetilde{u}(x,t)=\lambda^{-2}u(\lambda^{-1}x,\lambda^{-3})$ solves the Cauchy problem for initial data $\tilde{u}_0, \tilde{v}_0$ and $\tilde{w}_0$ on time interval $0\leq t\leq \lambda^3$.

\textbf{Step 5. Proof of locally Lipschitz continuity of map data-to-solution.} 

Let $\{(u_{0n},v_{0n},w_{0,n})\}_{n\in \{1,2\}}$ two initial data in $H^s(\mathcal{Y})$ such that $\|u_{0n}\|+\|v_{0n}\|+\|w_{0n}\|< \delta$, $(i=1,2)$ where $\delta$ is sufficiently small.

Let $(u_n,v_n,w_n)$ $(n=1,2)$ the solution of Cauchy problem \eqref{KDV}-\eqref{initial} with 1 boundary condition on the space $C([0,1]:H^s(\mathcal{Y}))$ with initial data $(u_{0n},v_{0n},w_{0n})$. According to Step 4 the lifespans of these solutions is $[0,1]$.

By using the arguments used in Step 4 we have that
\begin{equation}
\begin{split}
\|(u_2,v_2,w_2)-(u_1,v_1,w_1)\|_{Z|_{[0,1]}}\leq c\|(u_{02},v_{02},w_{02})-(u_{01},v_{01},w_{01})\|_{H^s(\mathcal{Y})}\\
+c(\|(u_2,v_2,w_2)+(u_1,v_1,w_1)\|_{Z|[0,1]})\|(u_2,v_2,w_2)-(u_1,v_1,w_1)\|_{Z|[0,1]},
\end{split}
\end{equation}
where $Z|_{[0,1]}$ denotes the restrictions of functions of $Z$ in the interval $[0,1].$
In a ball of $Z|_{[0,1]}$ we have that
\begin{equation}
\|(u_2,v_2,w_2)-(u_1,v_1,w_1)\|_{Z|_{[0,1]}}\leq c\|(u_{01},v_{02},w_{02})-(u_{01},v_{01},w_{01})\|_{H^s(\mathcal{Y})}.
\end{equation}
which completes the proof for small data assumptions. The local Lipschitz continuity for any data can be showed by a scaling argument.

\section{Proof of Corollary \ref{theorem}}\label{section7}

 By using Theorem \ref{theorem0} given a regularity index $s$ it suffices to get scalars $\lambda_i(s)$ satisfying \eqref{sclarsconditions} such that the matrix $M(\lambda_1,\lambda_2,\lambda_3,\lambda_4)$ given by \eqref{matrix} is invertible. These choices of scalars is a crucial point to get Corollary \ref{theorem}. We will divide this analysis in 2 cases.

\textbf{Case 1. Regularity: $\mathbf{-\frac{1}{2}<s<1,\ \text{for}\  s\neq\frac12}$.}

Taking $\lambda_i=0$ for $i=1,3,4$ and $0<\lambda_2=\frac{3}{\pi}\epsilon<<1$, a simple computations gives that the determinant of $\mathbf{M}$ is given by
\begin{equation}
\text{det}\, \mathbf{M}\left(0,\frac{3}{\pi}\epsilon,0,0\right)=2\sqrt{3}\ \alpha_2\ \alpha_3\ \text{sin}(\epsilon)\left( 1+\frac{1}{\alpha_2^2}+\frac{1}{\alpha_3^2}+\frac{\beta_3}{\alpha_3}+\frac{\beta_2}{\alpha_2}\right)\neq0,
\end{equation}
where we have used the hypothesis of Corollary \ref{theorem} about the parameters $\alpha_i$ and $\beta_i$ ($i\in\{1,2\}$), and the fact $0<\epsilon<<1$. Note that the condition \eqref{sclarsconditions} given in Theorem \ref{theorem0} is not valid for $\lambda=0$. Then, by a continuity argument we will take the parameters $\lambda_i\ (i=1,3,4)$ close to zero. In fact, for fixed $\alpha_2, \alpha_3, \beta_3$ and $\beta_4$ satisfying the hypothesis, we have that the function $\lambda \mapsto \text{det}\, \mathbf{M}(\lambda,\frac{3}{\pi}\epsilon,\lambda,\lambda)$ is continuous from $\R$ to $\mathbb{C}$. It follows that there exists a positive number $\delta(\epsilon)<<1$, depending of $\epsilon$, such that $\text{det}\, \mathbf{M}(\lambda,\frac{3\epsilon}{\pi},\lambda,\lambda)\neq0$, for $0<\lambda<\delta$. 

Thus, given $-\frac{1}{2}<s<1$ we can choice $(\lambda_1,\lambda_2,\lambda_3,\lambda_4)=(\lambda,\frac{3}{\pi}\epsilon,\lambda,\lambda)$ satisfying
\begin{equation}
\begin{split}
&0<\frac{3\epsilon}{2\pi}<\min \{s+\frac{1}{2},\frac12\}\\
&0<\lambda<\min\left\{\delta(\epsilon),s+\frac{1}{2}\right\}.
\end{split}
\end{equation}
Note that with this choice the all hypothesis of Theorem \ref{theorem0} part (i) are valid and the proof of Corollary \ref{theorem} on the \textbf{Case 1} is complete. 

\medskip

\textbf{Case 2. Regularity:} $\mathbf{1\leq s<\frac{3}{2}.}$ 

Taking $\lambda_i=\frac12$ for $i=1,3,4$ and $0<\lambda_2=\frac{1}{2}-\frac{3\epsilon}{\pi}$, for $0<\epsilon<<1$. A simple calculation shows that the determinant is given by
\begin{equation}
\text{det}\, \mathbf{M}\left(\frac12,\frac12-\frac{3\epsilon}{\pi},\frac12,\frac12\right)=2\sqrt{3}\ \alpha_2\ \alpha_3\  \text{sin}(\epsilon)\left( 1+\frac{1}{\alpha_2^2}+\frac{1}{\alpha_3^2}+\frac{\beta_3}{\alpha_3}+\frac{\beta_2}{\alpha_2}\right)\neq0,
\end{equation}
where we have used the hypothesis  $ \frac{1}{\alpha_2^2}+\frac{1}{\alpha_3^2}+\frac{\beta_3}{\alpha_3}+\frac{\beta_2}{\alpha_2}\neq-1 $ and the fact $0<\epsilon<<1$.

As the estimate condition \eqref{sclarsconditions} in part (i) of Theorem \eqref{theorem0} is not valid for $\lambda=\frac{1}{2}$, then we shall make a few perturbation in $\lambda$.

For fixed $\alpha_2,\ \alpha_3,\ \beta_3$ and $\beta_4$ satisfying the hypothesis, we have that the function $\lambda \mapsto \text{det}\, \mathbf{M}(\lambda,1-\frac{3}{\pi}\epsilon,\lambda,\lambda)$ is continuous from $\R$ to $\mathbb{C}$. It follows that there exists a positive number $\delta(\epsilon)<<1$, depending of $\epsilon$, such that $\text{det}\, \mathbf{M}(\lambda,\frac12-\frac{3}{2}\epsilon,\lambda,\lambda)\neq0$, for $0<\frac12-\lambda<\delta$. 

Thus, given $1<s<\frac{3}{2}$ we can choice $(\lambda_1,\lambda_2,\lambda_3,\lambda_4)=(\lambda,\frac12-\frac{3}{\pi}\epsilon,\lambda,\lambda)$ satisfying
\begin{equation}
\begin{split}
&s-1<\frac12-\frac{3\epsilon}{\pi}<s+\frac{1}{2},\\
&\max\left\{s-1,\frac12-\delta\right\}<\lambda<s+\frac{1}{2}.
\end{split}
\end{equation}
This finish the proof of Corollary \ref{theorem}.

\section{Proof of Corollary \ref{theorem1}}\label{section8}
For the regularity $-\frac12<s<1$ the result follows a similar idea of the proof of Corollary \ref{theorem}, by using the fact
\begin{equation}
\text{det}\, \mathbf{M}\left(0,\frac{3}{\pi}\epsilon,0,0\right)=2\sqrt{3}\alpha_2\alpha_3 \text{sin}(\epsilon)\left( 1+\frac{1}{\alpha_2^2}+\frac{1}{\alpha_3^2}+\frac{\beta_3}{\alpha_3}+\frac{\beta_2}{\alpha_2}\right)\neq0,
\end{equation}

Similarly, the case $1\leq s<\frac32$ with $s\neq\frac12$ follows from the fact
\begin{equation}
\text{det}\, \mathbf{M}\left(\frac12,\frac12-\frac{3\epsilon}{\pi},\frac12,\frac12\right)=2\sqrt{3}\alpha_2\alpha_3 \text{sin}(\epsilon)\left( 1+\frac{1}{\alpha_2^2}+\frac{1}{\alpha_3^2}+\frac{\beta_3}{\alpha_3}+\frac{\beta_2}{\alpha_2}\right)\neq0.
\end{equation}

\section*{Acknowledgments}
The author wishes to thank the Centro de Modelamiento Matem\'atico (CMM) and Universidad de Chile and N\'ucleo Milenio CAPDE, for the financial support and nice scientific infrastructure that allowed to conclude the paper during his postdoctoral stay.

\end{document}